\documentclass[11pt,reqno]{amsart}
\usepackage[T1]{fontenc}
\usepackage{lmodern,amsmath,amssymb,amsthm,mathtools,mathrsfs}
\usepackage{microtype,booktabs,array,enumitem}
\usepackage[a4paper,top=27mm,bottom=27mm,left=27mm,right=27mm]{geometry}
\usepackage[hidelinks]{hyperref}
\usepackage{bookmark}
\hypersetup{pdftitle={Quasipositive and almost positive curvature on finite quotients of Wallach spaces},pdfauthor={},pdfsubject={Explicit metrics, complete flat-plane descriptions, and finite-quotient topology in dimensions 6, 12 and 24}}
\newtheorem{maintheorem}{Theorem}

\newtheorem{theorem}{Theorem}[section]
\newtheorem{proposition}[theorem]{Proposition}
\newtheorem{lemma}[theorem]{Lemma}
\newtheorem{corollary}[theorem]{Corollary}
\theoremstyle{definition}
\theoremstyle{remark}\newtheorem{remark}[theorem]{Remark}
\numberwithin{equation}{section}
\DeclareMathOperator{\SU}{SU}\DeclareMathOperator{\Sp}{Sp}
\DeclareMathOperator{\SO}{SO}\DeclareMathOperator{\UU}{U}
\DeclareMathOperator{\Spin}{Spin}\DeclareMathOperator{\Aut}{Aut}
\DeclareMathOperator{\Der}{Der}\DeclareMathOperator{\Ad}{Ad}
\DeclareMathOperator{\Lie}{Lie}\DeclareMathOperator{\diag}{diag}
\DeclareMathOperator{\tr}{tr}\DeclareMathOperator{\Real}{Re}
\DeclareMathOperator{\Imag}{Im}\DeclareMathOperator{\rank}{rank}
\DeclareMathOperator{\Span}{span}\DeclareMathOperator{\secc}{sec}
\DeclareMathOperator{\Vol}{Vol}
\DeclareMathOperator{\Gr}{Gr}
\DeclareMathOperator{\im}{im}
\DeclareMathOperator{\Stab}{Stab}
\newcommand{\R}{\mathbb R}\newcommand{\C}{\mathbb C}
\newcommand{\HH}{\mathbb H}\newcommand{\OO}{\mathbb O}
\newcommand{\F}{\mathbb F}
\newcommand{\Qrat}{\mathbb Q}\newcommand{\CP}{\mathbb{CP}}
\newcommand{\OP}{\mathbb{OP}}
\newcommand{\RP}{\mathbb{RP}}
\newcommand{\gfr}{\mathfrak g}\newcommand{\kfr}{\mathfrak k}
\newcommand{\lfr}{\mathfrak l}\newcommand{\pfr}{\mathfrak p}
\newcommand{\mfr}{\mathfrak m}\newcommand{\afr}{\mathfrak a}
\newcommand{\qfr}{\mathfrak q}\newcommand{\tfr}{\mathfrak t}
\newcommand{\sfr}{\mathfrak s}\newcommand{\frso}{\mathfrak{so}}
\newcommand{\ZZ}{\mathcal Z}\newcommand{\BB}{\mathscr B}
\newcommand{\CC}{\mathscr C}\newcommand{\PP}{\mathscr P}
\newcommand{\NN}{\mathcal N}\newcommand{\calD}{\mathcal D}
\newcommand{\calE}{\mathcal E}\newcommand{\calN}{\mathcal N}
\newcommand{\ii}{\mathbf i}\newcommand{\jj}{\mathbf j}\newcommand{\kk}{\mathbf k}
\newcommand{\ip}[2]{\langle #1,#2\rangle}

\newcommand{\doi}[1]{\href{https://doi.org/#1}{\nolinkurl{#1}}}
\newcommand{\arxiv}[1]{\href{https://arxiv.org/abs/#1}{arXiv:\nolinkurl{#1}}}
\setlist[enumerate]{label=\textup{(\roman*)},itemsep=3pt,topsep=5pt}
\allowdisplaybreaks[2]
\title[Finite quotients of Wallach spaces]{Quasipositive and almost positive sectional curvature\\ on finite quotients of Wallach spaces}
\author{Minbo Gao$^1$}
\author{Yuhang Liu$^{*2}$}
\email{yuhang.liu02@xjtlu.edu.cn}
\date{November 9, 2026}
\subjclass[2020]{53C20, 53C30, 57S15}
\keywords{Almost positive sectional curvature, quasipositive curvature, Wallach spaces, Cheeger deformation, normal biquotients, triality, fundamental group}

\begin{document}
\footnotetext[1]{Institute of Software, Chinese Academy of Sciences}
\footnotetext[2]{Department of Applied Mathematics, Xi'an Jiaotong-Liverpool University}
\begin{abstract}
We construct nonnegatively curved metrics on the cyclic and full permutation quotients of the Wallach flag manifolds in dimensions $6$, $12$, and $24$. The metrics are quasipositively curved in all three dimensions and almost positively curved in dimensions $12$ and $24$. The cyclic quotients are orientable with fundamental group $C_3$, and their further full permutation quotients have fundamental group $S_3$; none admits any metric of everywhere positive sectional curvature. For the almost-positive metrics we describe every zero-curvature plane and prove that the zero-point loci have Hausdorff dimensions exactly $11$ and $23$, respectively. A one-step deformation has, in each dimension, an explicit zero locus with nonempty interior. Further nested deformations reduce it to a measure-zero set while leaving a common codimension-one stratum. The argument combines normal forms for horizontal commuting planes with quaternionic rank conditions and a dimension estimate for commuting compressions of skew-symmetric operators. We also determine the possible degrees of arbitrary free finite quotients of the three flag covers and derive order bounds for fundamental groups when the Euler characteristic is nonzero. We have used ChatGPT-6-astra for mathematical exploration and verified and rewritten all the results. 
\end{abstract}
\maketitle
\tableofcontents

\section{Introduction}\label{sec:introduction}
A metric has \emph{quasipositive sectional curvature} if its sectional curvature is nonnegative everywhere and all two-planes at one point have strictly positive curvature. It has \emph{almost positive sectional curvature} if its sectional curvature is nonnegative everywhere and all two-planes at every point of an open dense subset have strictly positive curvature. Throughout this paper the strict inequalities concern all planes at a common base point. For a nonnegatively curved metric $g$ on a closed manifold $M$, write
\begin{equation}\label{eq:Zdefinition}
 \ZZ(g)=\{x\in M:\text{some }\sigma\in\Gr_2(T_xM)\text{ satisfies }\secc_g(\sigma)=0\}.
\end{equation}
This is the \emph{zero-point locus}, as distinguished from the zero-plane locus in $\Gr_2(TM)$. It is closed. Almost positivity is equivalent to its having empty interior. In the examples constructed below it has, more strongly, codimension one in Hausdorff dimension and zero Riemannian volume.

\textbf{Acknowledgement:} the authors thank Xi’an Jiaotong-Liverpool University for providing financial and academic support. We also thank Professor Wolfgang Ziller and Professor Martin Kerin for helpful discussions. 

\subsection{Background and the role of finite quotients}
The difference between these curvature conditions is geometric and topological. A homogeneous metric that is quasipositively curved is positively curved everywhere, because isometries carry a positive point to every point. A homogeneous \emph{manifold}, however, may carry an inhomogeneous metric with quasipositive or almost positive curvature even when its topology excludes positive curvature. Finite quotients are particularly useful in making this distinction explicit.

The classical construction of Gromoll and Meyer gave an exotic seven-sphere with nonnegative sectional curvature and a point of positive sectional curvature~\cite{GromollMeyer}. Petersen and Wilhelm subsequently constructed almost-positive metrics on the unit tangent bundle of $S^4$ and a six-dimensional quotient~\cite{PetersenWilhelm}. Almost positivity on the Gromoll--Meyer sphere was established by Wilhelm~\cite{Wilhelm}; Eschenburg and Kerin later gave a simpler construction and proof~\cite{EschenburgKerin}. These examples demonstrated that the set of flat planes could be reduced substantially without preserving homogeneity.

Wilking's work~\cite{Wilking} introduced a broad normal-biquotient construction of almost-positive metrics. His families include projectivized tangent bundles over real, complex, and quaternionic projective spaces and certain lens-space bundles. Some of the real projective tangent bundles are odd-dimensional and nonorientable, so Synge's theorem excludes positive sectional curvature on their underlying manifolds. Thus almost positivity was already known to be strictly more flexible than positive curvature at the level of smooth manifolds. Tapp reformulated these constructions in the language of homogeneous bundles, gave effective quasipositivity conditions, and related the relevant obstruction to the fatness of a bundle~\cite{Tapp}. Kerr and Tapp extended this approach through positive triples and obtained additional homogeneous bundles and finite cyclic quotients of $T^1S^6$~\cite{KerrTapp}.

Later developments sharpened the analysis of flat planes in several directions. Kerin constructed almost-positive metrics on nontrivial $S^7$-bundles over $S^4$ and $\CP^3$~\cite{Kerin}. DeVito and Nance proved almost positivity on the irreducible rank-two symmetric space $G_2/\UU(2)\cong\Gr_2^+(\R^7)$~\cite{DeVitoNance}. DeVito and Johnson's analysis of Eschenburg spaces distinguishes the curvature properties of different natural inducing metrics and illustrates why quasipositivity of an analytic quotient metric does not imply almost positivity~\cite{DeVitoJohnson}. More recently, DeVito and West constructed infinitely many almost-positive generalized Eschenburg spaces in each dimension $4n-1$, $n\ge3$; for $n\ge6$, infinitely many are not homotopy equivalent to homogeneous spaces~\cite{DeVitoWest}. The last result is cited as a preprint. The present paper addresses a different feature: almost positivity on even-dimensional finite quotients whose fundamental groups rule out every positively curved metric.

Our manifolds come from the classical Wallach spaces. For $d=2,4,8$, respectively, let $\F=\C,\HH,\OO$ and use the following subgroup chains:
\begin{equation}\label{eq:groupstable}
\begin{array}{c|c|c|c}
 d&G_d&K_d&L_d\\ \hline
 2&\SU(3)&T^2&S(\UU(2)\times\UU(1))\\
 4&\Sp(3)&\Sp(1)^3&\Sp(2)\times\Sp(1)\\
 8&F_4&\Spin(8)&\Spin(9).
\end{array}
\end{equation}
Put $W_d=G_d/K_d$. The homogeneous fibration
\begin{equation}\label{eq:Wallachbundle}
 S^d=L_d/K_d\longrightarrow W_d\longrightarrow G_d/L_d=\F P^2
\end{equation}
has total dimension $3d$. Wallach proved that these total spaces admit homogeneous metrics with positive sectional curvature~\cite{Wallach}. They are also regular leaves of Cartan's three-principal-curvature isoparametric families. Their cyclic quotients occur as Gauss images in the work of Ma and Ohnita~\cite[Section 5]{MaOhnita}. In particular, the smooth manifolds studied here, and normal homogeneous metrics of nonnegative curvature on them, predate the present constructions. Their occurrence in isoparametric geometry does not identify their spherical induced metrics with the positive Wallach metrics.

The normalizer of $K_d$ has component group $S_3$, acting by permutation of an ordered projective frame. Define
\begin{equation}\label{eq:XY}
 X_d=G_d/(K_d\rtimes C_3),\qquad
 Y_d=G_d/(K_d\rtimes S_3).
\end{equation}
The action is on the right. The cyclic action preserves orientation, whereas a transposition reverses orientation. Thus $X_d$ is orientable with fundamental group $C_3$, and $Y_d$ is nonorientable with fundamental group $S_3$. Synge's theorem~\cite{Synge} excludes positive sectional curvature on both. In particular, no positively curved Wallach metric is invariant under the cyclic action. A metric suitable for descent must be constructed separately.

\subsection{Main results}
Fix an $\Ad(G_d)$-invariant inner product $Q$ on $\mathfrak g_d$. We use $Q(X,Y)=-\tfrac12\Real\tr(XY)$ for the matrix groups and the negative Killing form for $F_4$. All group metrics below are extended by \emph{right} translations. This convention makes the whole disconnected subgroup in \eqref{eq:XY} act isometrically.

\begin{maintheorem}\label{thm:existence}
For each $d\in\{2,4,8\}$, both $X_d$ and $Y_d$ admit explicit quasipositively curved metrics. For $d=4$ and $d=8$, they admit explicit almost-positively curved metrics. None of these six manifolds admits a metric with everywhere positive sectional curvature. Moreover,
\[
 \pi_1(X_d)=C_3,\quad\pi_1(Y_d)=S_3,\qquad
 H^*(X_d;\Qrat)\cong H^*(S^{3d};\Qrat),\quad
 H^*(Y_d;\Qrat)\cong H^*(\mathrm{pt};\Qrat).
\]
The map $X_d\to Y_d$ is the orientation double cover.
\end{maintheorem}

The metrics have simple inverse metric operators. For the initial one-step deformation let $P_L$ be orthogonal projection to $\mathfrak l_d$, and put
\begin{equation}\label{eq:initialmetric}
 \Phi_t^{-1}=I+(t^{-1}-1)P_L,\qquad 0<t<1.
\end{equation}
It induces a metric $g_t$ on every finite quotient under consideration. For the stronger quaternionic construction use
\begin{equation}\label{eq:quatchain}
 \Sp(3)\supset L_1=\Sp(2)\Sp(1)
 \supset L_2=\Sp(1)^3\supset L_3=\Delta\Sp(1),
 \qquad
 \Phi_{\HH}^{-1}=I+s_1P_1+s_2P_2+s_3P_3,
\end{equation}
where $s_1,s_2,s_3>0$. For the exceptional construction use the coordinate chain
\begin{equation}\label{eq:exceptionalchain}
 F_4\supset\Spin(9)\supset\Spin(8)\supset\cdots\supset\Spin(3),
 \qquad
 \Phi_{\OO}^{-1}=I+\sum_{j=3}^{9}s_jP_j,
 \quad s_j>0.
\end{equation}
The embeddings in \eqref{eq:exceptionalchain} are specified in Section~\ref{sec:exceptional}. In each case the metric at the group identity is $Q(\Phi\cdot,\cdot)$.

To describe the flat regions uniformly, let $r_1,r_2,r_3$ be the diagonal coordinates of the projective frame relative to a fixed point. For $\F=\C,\HH$ these are $r_i(aK_d)=|a_{3i}|^2$. In the exceptional case, with $E_i$ the coordinate idempotents of the Albert algebra, they are $r_i(aK_d)=\langle E_3,aE_i\rangle$. In all cases $r_i\ge0$, $\sum r_i=1$, and
\begin{equation}\label{eq:Psi}
 \Psi=\sum_{i=1}^3r_i^2
\end{equation}
is well defined on both finite quotients.

\begin{maintheorem}\label{thm:loci}
For every $d\in\{2,4,8\}$, the one-step metric satisfies
\begin{equation}\label{eq:oldall}
 \ZZ(g_t)=\{\Psi\le\tfrac12\}.
\end{equation}
This set has nonempty interior. The metrics defined by \eqref{eq:quatchain} and \eqref{eq:exceptionalchain}, denoted $g_{\HH}$ and $g_{\OO}$, satisfy
\begin{equation}\label{eq:mainloci}
 \{\Psi=\tfrac12\}\subset\ZZ(g_{\F})\subset\{\Psi\le\tfrac12\},
 \qquad
 \dim_{\mathrm H}\ZZ(g_{\F})=3d-1,
 \qquad \Vol_{g_{\F}}\ZZ(g_{\F})=0,
\end{equation}
for $(\F,d)=(\HH,4),(\OO,8)$. Their zero-point loci and normal horizontal flat-plane data are independent of the positive strengths $s_j$. A complete parametrization in the quaternionic case is given in Theorem~\ref{q:thm:exact} and Proposition~\ref{q:prop:allplanes}; the exceptional case is described in Proposition~\ref{e:prop:exactlocus}.
\end{maintheorem}

The inequalities in \eqref{eq:mainloci} are not meant to replace those complete descriptions. In particular, the zero locus need not be just the level set $\Psi=1/2$. The almost-positive metrics are different from the one-step quasipositive metrics: the additional deformations remove an open region of flat points, but a codimension-one stratum necessarily survives within these constructions.

The topology also gives sharp restrictions that do not assume any symmetry of the metric or the action.
\begin{maintheorem}\label{thm:degrees}
If a finite group $\Gamma$ acts freely and smoothly on $W_d$, where $d=2,4,8$, then
\[
 |\Gamma|\in\{1,2,3,6\}.
\]
If the quotient is orientable, then $|\Gamma|\in\{1,3\}$. The maximal orders $3$ and $6$ are attained by the orientable and nonorientable examples in Theorem~\ref{thm:existence}. A nontrivial further finite free quotient of $X_d$ has degree two and is nonorientable.
\end{maintheorem}

For general quasipositively curved closed $n$-manifolds, the fundamental group is finite by the Cheeger--Gromoll splitting theorem~\cite{CheegerGromoll}. Combining this with Gromov's Betti-number bound~\cite{Gromov} gives
\[
 |\pi_1(M)|\,|\chi(M)|\le B_n,
\]
where $B_n$ depends only on $n$. Section~\ref{sec:orders} derives the precise consequences, including bounds in dimension four and for six-dimensional biquotient universal covers. These statements do not claim a dimension-only group-order bound when the Euler characteristic vanishes.

\subsection{Proof strategy and scope}
The common construction is a quotient of a product of compact groups with bi-invariant metrics. The submersion formula gives nonnegative curvature. Tapp's flat-plane theorem~\cite{TappFlats} then makes the flatness test exact: a quotient plane is flat precisely when its normal horizontal lift commutes in the group and after every projection occurring in the nested deformation. This criterion, established in the form needed here in Section~\ref{sec:prelim}, is used throughout.

The main algebraic step is to classify all normal horizontal commuting two-planes. In the complex, quaternionic, and exceptional cases, they form a single isotropy orbit of a Fourier-conjugate two-torus. This reduces the first deformation to a row or torus-isotropy calculation, yielding \eqref{eq:oldall}. In dimension twelve the two further bracket conditions become a rank condition on three imaginary quaternionic axes. In dimension twenty-four the successive spin-group projections become commuting compressions of two skew-symmetric operators of ranks two and six. Their admissible flags have dimension at most $29$. After subtracting stabilizers and accounting for base-point parameters, this gives the sharp bound $23$ for the exceptional zero locus. The rank-drop branches are included in both arguments.

There is a genuine difference in dimension six. The complex diagonal algebra is abelian, and the additional mechanism used over $\HH$ has no complex analogue. Theorem~\ref{thm:complexrestriction} proves that every finite nested connected-subgroup deformation of a bi-invariant metric on $\SU(3)$, followed by the indicated right quotient, has an open set of flat points. Proposition~\ref{prop:conformal} also shows that a nonnegatively curved conformal change of any normal-biquotient metric retains every old zero plane. Neither statement rules out arbitrary almost-positive metrics on $X_2$ or $Y_2$.

The constructions strengthen the known curvature on the Gauss-image quotients, not their classification as smooth manifolds. The previous literature supplies the Wallach spaces, their isoparametric interpretation, the submersion machinery, and the topological obstructions. The contributions proved here are the equivariant curvature constructions, the complete normal-form and zero-locus analyses, and the dimension estimates that establish almost positivity on the finite quotients.

\section{Preliminaries: submersions, nested deformations, and dimension}\label{sec:prelim}
We use $R(X,Y)Z=\nabla_X\nabla_YZ-\nabla_Y\nabla_XZ-\nabla_{[X,Y]}Z$ and $\secc(X\wedge Y)=\langle R(X,Y)Y,X\rangle/|X\wedge Y|^2$. All groups are compact, and a closed manifold is compact without boundary. Orthogonal complements and projections in Lie algebras are taken with respect to $Q$.

\subsection{Curvature facts used throughout}
For a bi-invariant group metric, the curvature numerator is $\tfrac14|[X,Y]|^2$. For a Riemannian submersion $\pi:E\to B$ and horizontal vectors $X,Y$, O'Neill's formula~\cite{ONeill} reads
\begin{equation}\label{eq:ONeill}
 \langle R^B(d\pi X,d\pi Y)d\pi Y,d\pi X\rangle
 =\langle R^E(X,Y)Y,X\rangle+3|\mathcal A_XY|^2.
\end{equation}
In particular a submersion preserves a nonnegative lower sectional-curvature bound. Necessity of commuting-lift equations follows at once from this identity, but sufficiency requires control of the last term. We use the following theorem instead of repeating that argument in each example.

\begin{theorem}[Tapp~{\cite[Theorem 1.1]{TappFlats}}]\label{thm:Tappflat}
Let a compact Lie group with a bi-invariant metric admit a Riemannian submersion onto $B$. Every horizontal zero-curvature plane projects to a zero-curvature plane of $B$. Every zero-curvature plane of $B$ exponentiates to a complete totally geodesic flat immersion of $\R^2$.
\end{theorem}
A product of compact groups with scaled bi-invariant metrics is again a compact Lie group with a bi-invariant metric, so the theorem applies to all total spaces below. No connectedness assumption on the fibers of the final quotient is required.

\subsection{An exact criterion for a nested deformation}
Let $G\supset L_1\supset\cdots\supset L_s$ be closed connected subgroups. Write $P_j$ for projection to $\mathfrak l_j$ and choose strengths $c_j>0$. Set
\begin{equation}\label{eq:generalPhi}
 \Phi^{-1}=I+\sum_{j=1}^s c_jP_j,
 \qquad \overline g_e(U,V)=Q(\Phi U,V),
\end{equation}
and extend $\overline g$ by right translations. This is the iterated Cheeger construction~\cite{Cheeger}; it can also be realized directly by
\begin{equation}\label{eq:productmap}
 (G,Q)\times\prod_{j=1}^s(L_j,c_j^{-1}Q)
 \longrightarrow (G,\overline g),\qquad
 (a,b_1,\ldots,b_s)\longmapsto b_s^{-1}\cdots b_1^{-1}a.
\end{equation}
The fibers are the orbits of
\[
 (h_1,\ldots,h_s)\cdot(a,b_1,\ldots,b_s)
 =(h_1a,h_1b_1h_2^{-1},\ldots,h_{s-1}b_{s-1}h_s^{-1},h_sb_s).
\]
The inclusions make this a free isometric action. At $(a,e,\ldots,e)$, right trivialization identifies the horizontal lift of $\Phi^{-1}V$ with
\begin{equation}\label{eq:lift}
 (V,-c_1P_1V,\ldots,-c_sP_sV).
\end{equation}
Indeed, orthogonality to the successive infinitesimal action fields gives these components; applying the differential of \eqref{eq:productmap} gives $\Phi^{-1}V$, and the squared norm is $Q(V,\Phi^{-1}V)$. Right translation on the first factor descends to right translation on $G$, proving the stated metric formula.

Let $H\subset G$ be a closed subgroup with identity component $K$ and Lie algebra $\mathfrak k$. Its right action is isometric. At $a\in G$, its vertical space in right trivialization is $\Ad_a\mathfrak k$. Thus \eqref{eq:lift} is horizontal for the composite map to $G/H$ exactly when $V\perp_Q\Ad_a\mathfrak k$.

\begin{proposition}[Exact flat-plane criterion]\label{prop:flatcriterion}
The quotient metric on $G/H$ is nonnegatively curved. Every zero plane at $aH$, and only such a plane, is obtained by projecting
\begin{equation}\label{eq:projectedplane}
 \Span\{dR_a(\Phi^{-1}V),dR_a(\Phi^{-1}W)\},
\end{equation}
where $V,W\in(\Ad_a\mathfrak k)^\perp$ are independent and
\begin{equation}\label{eq:brackets}
 [V,W]=0,\qquad [P_jV,P_jW]=0\quad(1\le j\le s).
\end{equation}
Consequently the normal horizontal zero planes, and their base-point locus, are independent of the positive strengths $c_j$.
\end{proposition}
\begin{proof}
Compose \eqref{eq:productmap} with the right quotient. The resulting map is a Riemannian submersion from a product of bi-invariant groups. Its horizontal lift is \eqref{eq:lift}. If a quotient plane is flat, \eqref{eq:ONeill} forces its lift to have zero curvature. The bi-invariant curvature formula in every factor gives exactly \eqref{eq:brackets}. Conversely those equations make the horizontal lifted plane flat, and Theorem~\ref{thm:Tappflat} proves that its projection is flat. All $c_j$ are positive, so their values do not affect the equations.
\end{proof}

The word ``normal'' in this criterion refers to $Q$-horizontality of $V,W$, not to horizontality of the unmodified vectors for $\overline g$. The actual horizontal plane is \eqref{eq:projectedplane}. This distinction will be maintained in the zero-plane parametrizations.

\subsection{A finite-extension positivity criterion}\label{subsec:fatness}
The following observation isolates the local positivity mechanism common to the three flag spaces. It is compatible with, but distinct from, the positive homogeneous-bundle criterion in Tapp~\cite[Theorem 2.2 and Section 3]{Tapp}.
\begin{proposition}\label{prop:finiteextension}
Let $K\subset L\subset G$ and write $\mathfrak g=\mathfrak k\oplus\mathfrak m\oplus\mathfrak p$ orthogonally, with $\mathfrak l=\mathfrak k\oplus\mathfrak m$. Suppose $(G,L)$ is a symmetric pair, commuting vectors in each of $\mathfrak m$ and $\mathfrak p$ are real-linearly dependent, and
\begin{equation}\label{eq:fatness}
 [A,P]\ne0\qquad(0\ne A\in\mathfrak m,\quad0\ne P\in\mathfrak p).
\end{equation}
For every closed $H\subset G$ with identity component $K$, the right-invariant one-step metric $tQ|_{\mathfrak l}+Q|_{\mathfrak p}$, $0<t<1$, induces strictly positive sectional curvature on every plane at $eH$ and nonnegative curvature everywhere.
\end{proposition}
\begin{proof}
Nonnegativity follows from Proposition~\ref{prop:flatcriterion}. At $eH$, a normal horizontal commuting pair can be written $V=A+P$, $W=B+R$, with $A,B\in\mathfrak m$ and $P,R\in\mathfrak p$. Flatness would give $[A,B]=0$ and $[V,W]=0$. The symmetric splitting gives $[P,R]=0$ and $[A,R]+[P,B]=0$. Hence $A=aA_0$, $B=bA_0$, $P=cP_0$, $R=dP_0$. If either component pair vanished, $V,W$ would be dependent. Otherwise independence gives $ad-bc\ne0$, whereas the mixed equation is $(ad-bc)[A_0,P_0]=0$, contradicting \eqref{eq:fatness}.
\end{proof}
For the chains \eqref{eq:groupstable}, the rank-one symmetric bases and sphere fibers give the two commuting-vector conditions; the mixed condition is the fatness of the Wallach sphere bundles~\cite{BerardBergery,Tapp,Wallach}. These known facts already give quasipositivity of all the finite quotients. The complete zero-locus calculations below provide an independent, more precise verification.

\subsection{Finite quotients and dimension estimates}
A finite isometric covering preserves curvature pointwise. In particular, if $p:\widetilde M\to M$ is such a covering, then
\begin{equation}\label{eq:coverZ}
 \ZZ(p^*g)=p^{-1}\ZZ(g).
\end{equation}
Local isometry and a finite cover by evenly covered coordinate charts show that Hausdorff dimensions of these two loci agree. Quasipositivity, almost positivity, and the assertion that the zero locus has volume zero are therefore equivalent upstairs and downstairs. We will perform calculations on the ordered cover when convenient.

We use the standard dimension rules for real semialgebraic sets~\cite[Section 2.8]{BCR}: finite unions take the maximum dimension; a family has dimension at most the parameter dimension plus the largest fiber dimension; and an image has dimension at most the domain dimension minus a uniform lower bound for its fibers. Such sets have finite smooth stratifications, and their dimension equals their Hausdorff dimension in smooth coordinate charts. All group representations, projection equations, and rank conditions used here are real algebraic. Trigonometric angle intervals are replaced by their semialgebraic sine--cosine arcs. Consequently these rules apply to the parametrizations, including singular branches. We do not infer a dimension bound merely by counting equations.

\section{Topology of the flag quotients}\label{sec:topology}
We first identify the finite actions and their topological obstructions. The discussion applies uniformly to $d=2,4,8$ and does not impose any restriction on other metrics the quotient might carry.

\subsection{The normalizer and its free action}
\begin{proposition}\label{prop:normalizers}
For each chain in \eqref{eq:groupstable},
\begin{equation}\label{eq:normalizerS3}
 N_{G_d}(K_d)=K_d\rtimes S_3.
\end{equation}
The induced right action of $S_3$ on $W_d$ is free.
\end{proposition}
\begin{proof}
For $\SU(3)$, a normalizer of the diagonal torus permutes its three joint eigenspaces, hence is monomial. The map $\sigma\mapsto\operatorname{sgn}(\sigma)P_\sigma$ is a homomorphism from $S_3$ into $\SU(3)$ and provides a splitting. For $\Sp(3)$, the diagonal subgroup $\Sp(1)^3$ has three inequivalent irreducible coordinate quaternionic modules; a normalizer permutes them. Quaternionic monomial matrices therefore give the full normalizer, and the real permutation matrices split it.

For $F_4=\Aut(J)$, use the standard Albert-algebra description of the frame stabilizer~\cite{Baez,Kollross}. The $K_d$-fixed subspace of $J$ is the diagonal algebra $\Span(E_1,E_2,E_3)$. Any normalizer preserves this algebra and permutes its primitive idempotents. The kernel of the permutation action is exactly $K_d$. Conversely real permutation matrices act on $J$ by $A\mapsto PAP^T$, preserve the Jordan product, and give all permutations. Their intersection with $K_d$ is trivial, proving \eqref{eq:normalizerS3}.

Finally, if $a nK_d=aK_d$ for $n\in N_{G_d}(K_d)$, then $n\in K_d$. Thus no nontrivial normalizer component fixes a point of $W_d$.
\end{proof}

\subsection{Universal covers, orientation, and cohomology}
The standard groups $G_d$ are simply connected, and $K_d$ is connected. The homotopy sequence of $K_d\to G_d\to W_d$ therefore shows that $W_d$ is simply connected. The fibration \eqref{eq:Wallachbundle} also yields its rational Poincar\'e polynomial. The base has rational cohomology in degrees $0,d,2d$, and the fiber in degrees $0,d$. Every differential in the Serre spectral sequence vanishes by degree considerations. Consequently
\begin{equation}\label{eq:flagBetti}
 P_{W_d}(t)=1+2t^d+2t^{2d}+t^{3d},\qquad \chi(W_d)=6.
\end{equation}
The fibration and its projective-plane interpretation are classical~\cite{Wallach,MaOhnita}; the use of the spectral sequence here determines only additive cohomology, which is all that is needed.

\begin{lemma}\label{lem:orientationcharacter}
The orientation character of the free $S_3$-action on $W_d$ is its sign character.
\end{lemma}
\begin{proof}
The connected isotropy gives $W_d$ an orientation. Consider the transposition interchanging the first two projective-frame entries and fixing the third. It fixes the base point of the fibration \eqref{eq:Wallachbundle} and acts on its fiber by exchanging the two complementary rank-one idempotents. On the sphere model of the fiber this is the antipodal map: subtracting half the identity identifies the two idempotents with opposite points of $S^d$. Its degree is $(-1)^{d+1}=-1$, since $d$ is even. The derivative on the base is the identity, so the transposition reverses the orientation of the total space. All transpositions are conjugate and generate $S_3$, proving the claim. In the exceptional case the same fiber model is the unit sphere in $\operatorname{Herm}_2(\OO)_0$.
\end{proof}

\begin{proposition}\label{prop:topologyXY}
The spaces $X_d$ and $Y_d$ are closed smooth manifolds with
\[
 \pi_1(X_d)\cong C_3,\qquad \pi_1(Y_d)\cong S_3,
 \qquad \chi(X_d)=2,\quad \chi(Y_d)=1.
\]
The first is orientable, the second nonorientable, and $X_d\to Y_d$ is the orientation double cover. Their rational cohomology is as stated in Theorem~\ref{thm:existence}.
\end{proposition}
\begin{proof}
Smoothness and the fundamental groups follow from Proposition~\ref{prop:normalizers} and simple connectivity of $W_d$. The orientation statement follows from Lemma~\ref{lem:orientationcharacter}. Euler characteristic is multiplicative under finite coverings, giving the displayed values.

Rational transfer identifies the cohomology of a finite free quotient with the invariant cohomology of its cover~\cite[Section 3.G]{Hatcher}; in particular it injects into that of the cover. Thus $X_d$ can have rational cohomology only in degrees $0,d,2d,3d$. Poincar\'e duality gives $b_d(X_d)=b_{2d}(X_d)$. Since all four degrees are even and $\chi(X_d)=2$, the middle Betti numbers vanish. This proves that $X_d$ is a rational homology sphere. Transfer for $X_d\to Y_d$ leaves only degrees zero and $3d$, and nonorientability eliminates the top rational cohomology of $Y_d$.
\end{proof}

\begin{corollary}\label{cor:noPositive}
None of the spaces $X_d$ or $Y_d$ admits a metric with everywhere positive sectional curvature.
\end{corollary}
\begin{proof}
By Synge's theorem~\cite{Synge}, an orientable even-dimensional closed positively curved manifold is simply connected. This contradicts $\pi_1(X_d)=C_3$. Any positive metric on $Y_d$ would lift to one on its orientation cover $X_d$, giving the same contradiction.
\end{proof}
This obstruction concerns arbitrary smooth metrics, not just invariant metrics. It also explains why the positive homogeneous metrics on $W_d$ cannot descend to $X_d$.

\subsection{Sharp restrictions on all finite quotients}
\begin{proof}[Proof of Theorem~\ref{thm:degrees}]
Let $\Gamma$ act freely and smoothly on $W_d$. Equation~\eqref{eq:flagBetti} and Euler multiplicativity give
\[
 6=|\Gamma|\,\chi(W_d/\Gamma).
\]
Thus $|\Gamma|$ divides six. If the quotient is orientable, rational transfer again confines its cohomology to degrees $0,d,2d,3d$. Poincar\'e duality gives
\[
 \chi(W_d/\Gamma)=2+2b_d(W_d/\Gamma),
\]
a positive even divisor of six. Therefore $|\Gamma|=1$ or $3$. The free $C_3$ and $S_3$ actions above realize the maximal orders. Their required curvature metrics follow from the subsequent sections.

If a finite group $\Lambda$ acts freely on $X_d$, then
$2=|\Lambda|\chi(X_d/\Lambda)$, so a nontrivial such action has order two. An orientable quotient would, by transfer from the rational homology sphere $X_d$, have Euler characteristic two; the quotient actually has Euler characteristic one. It is therefore nonorientable.
\end{proof}
In particular, the order bounds do not arise from a limitation of monomial matrix representatives. They apply to every smooth free action on these universal covers.

\section{Normal flats in the complex and quaternionic flags}\label{sec:matrixflats}
In this section $\F=\C$ or $\HH$, and $(G,K,L)$ is the corresponding chain in \eqref{eq:groupstable}; we also write $L_1=L$. Let $H=K\rtimes\Gamma$, where $\Gamma$ is $C_3$ or $S_3$. The calculations also apply to $H=K$. Matrices act on right quaternionic vector spaces in the quaternionic case, and $A^*$ denotes conjugate transpose. The normal horizontal space $\mathfrak k^\perp$ consists of the zero-diagonal skew-Hermitian matrices.

Let
\begin{equation}\label{q:eq:Fourier}
 \omega=e^{2\pi i/3},\qquad
 F=\frac1{\sqrt3}\begin{pmatrix}1&1&1\\1&\omega&\omega^2\\1&\omega^2&\omega\end{pmatrix},
 \qquad
 \afr=\{i\diag(\lambda_1,\lambda_2,\lambda_3):\lambda\in\R^3,\ \sum\lambda_j=0\}.
\end{equation}
In the quaternionic case fix $\C=\R+\R\ii\subset\HH$ and interpret $i=\ii$. The matrix $F$ is unitary, belongs to $\Sp(3)$ in the quaternionic case, and acts by conjugation on $\mathfrak{su}(3)$ in the complex case. Multiplying $F$ by a scalar makes its determinant one when needed, without changing its conjugation action.

The plane $\qfr_0=\Ad_F\afr$ consists of zero-diagonal matrices and is spanned over $\R$ by
\begin{equation}\label{q:eq:A0B0}
 A_0=\begin{pmatrix}0&1&-1\\-1&0&1\\1&-1&0\end{pmatrix},\qquad
 B_0=\begin{pmatrix}0&i&i\\i&0&i\\i&i&0\end{pmatrix}.
\end{equation}
Direct multiplication gives $[A_0,B_0]=0$.

\begin{lemma}[Rank-one off-block algebra]\label{q:lem:rankone}
For $v,w\in\F^2$, put
\[
 P(v)=\begin{pmatrix}0_{2\times2}&v\\-v^*&0\end{pmatrix}.
\]
If $[P(v),P(w)]=0$ and $v\neq0$, then $w$ is a real multiple of $v$.
\end{lemma}
\begin{proof}
The top-left block of the commutator is $wv^*-vw^*$. Its vanishing, applied to $v$, gives $w=vq$, where $q=(w^*v)/|v|^2$. Substitution gives $v(q-\bar q)v^*=0$, hence $q=\bar q\in\R$ after multiplication by $v^*$ and $v$.
\end{proof}

\begin{lemma}[Classification of commuting horizontal planes]\label{q:lem:normalclassification}
Over $\C$ and over $\HH$, every two-dimensional real abelian subspace of $\kfr^\perp$ has the form
\begin{equation}\label{q:eq:normalflat}
 \qfr=\Ad_D\qfr_0,\qquad D\in K.
\end{equation}
\end{lemma}
\begin{proof}
Matrices in $\kfr^\perp$ are precisely the zero-diagonal skew-Hermitian matrices. Write independent commuting elements as
\[
 X=\begin{pmatrix}0&x&y\\-\bar x&0&z\\-\bar y&-\bar z&0\end{pmatrix},\quad
 Y=\begin{pmatrix}0&u&v\\-\bar u&0&w\\-\bar v&-\bar w&0\end{pmatrix}.
\]
None of the three edge projections can vanish identically on the plane: otherwise the whole plane would lie in an off-block space from Lemma~\ref{q:lem:rankone}, contrary to independence. We may therefore choose $X$ with all three edge entries nonzero, avoiding three proper real linear subspaces. Diagonal unitary conjugation makes $x,z$ positive real numbers. In the complex case a common scalar corrects the determinant of this conjugating matrix and does not affect the conjugation.

The off-diagonal equations in $[X,Y]=0$ are
\begin{equation}\label{q:eq:edgeeqs}
 w=\frac zx u,\qquad v=\frac1x\bar u\,y,\qquad [y,u]=0.
\end{equation}
If $u$ were real, these equations would make $Y=(u/x)X$. Thus $\Imag u\neq0$. The first two diagonal equations give
\[
 (x-|y|^2/x)(u-\bar u)=0,\qquad
 (z^2/x-x)(u-\bar u)=0.
\]
Consequently $z=x=|y|$. Rescale $X$, subtract a real multiple of it from $Y$, and rescale $Y$ to obtain
\[
 x=z=1,\quad y=q,\quad |q|=1,\qquad
 u=I,\quad v=-Iq,\quad w=I,
 \qquad I^2=-1,\quad [q,I]=0.
\]
A common diagonal quaternionic conjugation makes $I=\ii$. In the complex case choose the sign of $Y$ to do the same. Hence $q=e^{i\varphi}$.

Replace the basis by
$X_\theta=\cos\theta X+\sin\theta Y$ and
$Y_\theta=-\sin\theta X+\cos\theta Y$.
Conjugating by $\diag(e^{-i\theta},1,e^{i\theta})$ makes the three edge entries of $X_\theta$ equal to
$1,q e^{-3i\theta},1$. Choose $\theta$ so that $q e^{-3i\theta}=-1$. The edge entries of $Y_\theta$ are then $i,i,i$. This is precisely \eqref{q:eq:A0B0}.
\end{proof}

\begin{corollary}\label{q:cor:eigenplane}
If $V,W\in(\Ad_a\kfr)^\perp$ are independent and commute, their span is
\begin{equation}\label{q:eq:eigenplane}
 \Span\{V,W\}=\Ad_B\afr,\qquad B=aDF,
 \quad\text{for some }D\in K.
\end{equation}
\end{corollary}

In the quaternionic case we will also need the normalizer of $\afr$:
\begin{equation}\label{q:eq:normalizer}
 \NN=\{EP,\ EP(\jj I_3):
 E=\diag(z_1,z_2,z_3),\ z_j\in S^1\subset\C,\ P\text{ a permutation matrix}\}.
\end{equation}
To check this, take a generic $\lambda$ with three distinct nonzero absolute values. Its quaternionic eigenspaces force a normalizing matrix to permute the three coordinate quaternionic lines. On each line it either preserves $i$ or sends it to $-i$. The relation $\lambda_1+\lambda_2+\lambda_3=0$ is the unique relation between the three coordinate linear forms on $\afr$; therefore all three signs must agree. This gives exactly \eqref{q:eq:normalizer}. Conversely each displayed matrix plainly normalizes $\afr$.

Thus replacing the eigenbasis $B$ in \eqref{q:eq:eigenplane} by another eigenbasis for the same plane multiplies it on the right by an element of $\NN$. Equivalently, if a canonical eigenbasis for that plane is $B_0$, its base point has the form
\begin{equation}\label{q:eq:basefromB}
 aH=B_0 nF^*H,\qquad n\in\NN.
\end{equation}
The right diagonal factor from \eqref{q:eq:eigenplane} has disappeared into $H$.

\subsection{The one-step zero locus}
For the single block deformation, Proposition~\ref{prop:flatcriterion} requires
\begin{equation}\label{q:eq:firstcondition}
 [V,W]=0,\qquad[V_{\lfr_1},W_{\lfr_1}]=0.
\end{equation}
Because $G/L_1$ is a symmetric pair, $[\pfr,\pfr]\subset\lfr_1$ and $[\lfr_1,\pfr]\subset\pfr$. If $[V,W]=0$, the second condition is equivalent to $[V_{\pfr},W_{\pfr}]=0$, and hence, by Lemma~\ref{q:lem:rankone}, to the $\pfr$-projection of their plane having rank at most one. Equivalently some nonzero element of their plane belongs to $\lfr_1$.

Let $B$ be the eigenbasis in \eqref{q:eq:eigenplane}, and write its last row as $b=(b_1,b_2,b_3)$. The last condition says that, for some nonzero real trace-zero $\lambda$, the matrix $B(i\diag\lambda)B^*$ preserves the line $\F e_3$. It is equivalent to
\begin{equation}\label{q:eq:roweigen}
 i\lambda_j\bar b_j=\bar b_j c\quad(j=1,2,3)
\end{equation}
for one scalar $c\in\Imag\F$.

\begin{lemma}[Last-row criterion]\label{q:lem:row}
For $\F=\C$, \eqref{q:eq:firstcondition} holds precisely when at least one of the entries $b_j$ is zero.
For $\F=\HH$, it holds precisely when either the row has only one nonzero entry, or it has exactly two nonzero entries $b_u,b_v$ and
\begin{equation}\label{q:eq:axisparallel}
 \frac{b_u i\bar b_u}{|b_u|^2}
 =\pm\frac{b_v i\bar b_v}{|b_v|^2}.
\end{equation}
\end{lemma}
\begin{proof}
In the complex case all nonzero row entries in \eqref{q:eq:roweigen} require the corresponding $\lambda_j$ to agree. Three nonzero entries would therefore force $\lambda=0$. With two or fewer entries, a nonzero trace-zero choice exists.

In the quaternionic case put $n_j=b_j i\bar b_j/|b_j|^2$ for nonzero entries. Equation \eqref{q:eq:roweigen} says $\lambda_j n_j=c$. If all three entries are nonzero, either all $\lambda_j=0$ or they have the same positive absolute value. Three signed numbers of that absolute value cannot sum to zero. With two nonzero entries, nonzero solutions exist exactly when their unit imaginary axes are parallel or antiparallel; with one entry they always exist.
\end{proof}

\begin{proof}[Proof of \eqref{eq:oldall} for $d=2,4$]
Set $s_j=|a_{3j}|$, so $s_1^2+s_2^2+s_3^2=1$. By Lemma~\ref{q:lem:row}, a zero plane requires a zero entry in the last row of $aDF$ for some diagonal $D\in K$. The triangle inequality then implies
\begin{equation}\label{q:eq:triangle}
 s_j\leq s_u+s_v\quad\text{for every }\{j,u,v\}=\{1,2,3\}.
\end{equation}
Conversely, if \eqref{q:eq:triangle} holds, choose complex numbers $c_j$ with $|c_j|=s_j$ and $c_1+c_2+c_3=0$. Choose unit diagonal entries $d_j$ so that $a_{3j}d_j=c_j$; a zero entry causes no difficulty. In the complex case multiply all $d_j$ by one common unit scalar to make the determinant one. The last row of $aDF$ is then complex and has a zero entry. Its remaining imaginary axes are equal, so Lemma~\ref{q:lem:row} gives a flat plane in both cases.

Finally, direct expansion, or Heron's identity, gives
\begin{align}\label{q:eq:Heron}
 &(s_1+s_2+s_3)(-s_1+s_2+s_3)
 (s_1-s_2+s_3)(s_1+s_2-s_3)\notag\\
 &\hspace{35mm}=2\sum_{u<v}s_u^2s_v^2-\sum_j s_j^4
 =1-2\sum_j s_j^4.
\end{align}
At most one of the last three factors can be negative. Thus \eqref{q:eq:triangle} is equivalent to $\Psi\leq1/2$. At the identity $\Psi=1$, whereas a representative with last-row norms $1/\sqrt3$ has $\Psi=1/3$. Continuity proves both quasipositivity and the existence of an open zero-curvature region.
\end{proof}

This proof also determines every zero plane of the one-step metric: take all $D\in K$ for which the last row of $B=aDF$ satisfies Lemma~\ref{q:lem:row}, conjugate the two-dimensional plane $\afr$ by $B$, and use \eqref{eq:projectedplane}. In particular \eqref{eq:oldall} is an equality, not an estimate for an exceptional region.

\subsection{A flat torus carrying the cyclic generator}\label{c:sec:flat}
The one-step six-dimensional metric has an embedded flat torus on which the cyclic fundamental group is represented by closed geodesics.

Let $T=K=T^2$, $H=T
times C_3$, and let $C$ be the cyclic permutation matrix. Set
\begin{equation}\label{c:eq:Fourier}
 F_0=\frac1{\sqrt3}\begin{pmatrix}1&1&1\\1&\omega&\omega^2\\1&\omega^2&\omega\end{pmatrix},
 \qquad a_0=e^{i\pi/6}F_0.
\end{equation}
One checks $F_0^*F_0=I$ and $\det F_0=-i$, so $a_0\in\SU(3)$. Every entry of $a_0$ has squared modulus $1/3$. Thus, for every $D=i\diag(d_1,d_2,d_3)\in\tfr$,
\begin{equation}\label{c:eq:zerodiagonal}
 (a_0^*Da_0)_{jj}=\frac i3(d_1+d_2+d_3)=0.
\end{equation}
Equivalently,
\begin{equation}\label{c:eq:torusorth}
 Q(\tfr,\Ad_{a_0}\tfr)=0.
\end{equation}

\begin{proposition}\label{c:prop:flat}
The torus $Ta_0\subset(G,\overline g)$ is horizontal and totally geodesic, with a flat induced metric. Its image in $(X_2,g_t)$ is an embedded totally geodesic flat two-torus. Moreover, this torus contains a closed geodesic representing a generator of $\pi_1(X_2)$.
\end{proposition}
\begin{proof}
First, $T$ is totally geodesic and flat in $(G,\overline g_t)$. For $D\in\tfr$, the Koszul formula for right-invariant vector fields gives
\[
 \overline Q_t(\nabla_{D^R}D^R,Z^R)
 =\overline Q_t([D,Z],D)=tQ([D,Z],D)=0.
\]
The last equality uses $\Ad$-invariance. Thus every one-parameter subgroup in $T$ is an ambient geodesic. Its induced metric is translation-invariant on the abelian group $T$, and hence flat. Right translation by $a_0$ preserves these properties.

In right trivialization, the vertical space of $G\to G/H$ at $a$ is $\Ad_a\tfr$. At a point $ta_0$ with $t\in T$, the tangent space to $Ta_0$ is $\tfr$ in the same trivialization. Because $\Ad_t$ fixes $\tfr$ pointwise,~\eqref{c:eq:torusorth} gives
\[
 \overline Q_t(\tfr,\Ad_{ta_0}\tfr)
 =tQ(\tfr,\Ad_{ta_0}\tfr)=0.
\]
This proves horizontality along the whole torus.

Horizontal ambient geodesics project to geodesics under a Riemannian submersion: this follows either from the horizontal connection identity in the submersion connection formula~\cite{ONeill} or directly from the first variation. The image of $Ta_0$ is therefore totally geodesic and locally isometric to a flat torus. More precisely, its stabilizer is $T\cap a_0Ha_0^{-1}$. Its Lie algebra is zero by~\eqref{c:eq:torusorth}, so the stabilizer is finite. The resulting torus quotient embeds as the compact $T$-orbit of $a_0H$.

Finally, direct multiplication gives
\[
 a_0Ca_0^{-1}=\diag(1,\omega,\omega^2)=:\tau\in T.
\]
Put $D_0=i\diag(0,2\pi/3,-2\pi/3)$, so $\exp(D_0)=\tau$, and consider
\begin{equation}\label{c:eq:nontrivialloop}
 \gamma(s)=\exp(sD_0)a_0H,\qquad 0\leq s\leq1.
\end{equation}
It is a geodesic because it is the projection of a horizontal geodesic in $Ta_0$. Since $\tau a_0=a_0C$, its endpoints agree in $X_2$. Its endpoint velocities agree as well: right multiplication by $C$ takes the initial lift and velocity to the final lift and velocity. Its lift to $G/T$ runs from $a_0T$ to $a_0CT$, so it represents a nontrivial deck transformation of order three. Thus it generates $\pi_1(X_2)$.
\end{proof}

In particular there are zero-curvature two-planes along this torus. The positive neighborhood of $o$ cannot meet it. This concretely illustrates how a quasipositive metric can have nontrivial loops supported in a region of zero sectional curvature.

\section{Almost positivity on the quaternionic quotients}\label{sec:quaternionic}
Write $g_*=g_{\HH}$ for the quaternionic metric. For the chain \eqref{eq:quatchain}, set
\[
 \pfr=\lfr_1^\perp,\quad
 \mfr=\lfr_1\ominus\lfr_2,\quad
 \afr_0=\lfr_2\ominus\lfr_3.
\]
Here and below $H=\Sp(1)^3\rtimes\Gamma$, with $\Gamma=C_3$ or $S_3$. The spaces have real dimensions $8,4,6$, while $\dim\lfr_3=3$. Formula \eqref{eq:quatchain} is equivalent to the particularly concrete metric
\begin{equation}\label{q:eq:weights}
 {\overline Q_*
 =Q|_{\pfr}+\tfrac12Q|_{\mfr}
 +\tfrac13Q|_{\afr_0}+\tfrac14Q|_{\Delta\mathfrak{sp}(1)}.}
\end{equation}
More generally, all three strengths may be chosen arbitrarily positive:
\begin{equation}\label{q:eq:family}
 \overline Q_{s_1,s_2,s_3}=Q|_{\pfr}
 +\frac{Q|_{\mfr}}{1+s_1}
 +\frac{Q|_{\afr_0}}{1+s_1+s_2}
 +\frac{Q|_{\lfr_3}}{1+s_1+s_2+s_3}.
\end{equation}
Proposition~\ref{prop:flatcriterion} shows that the zero-point locus does not depend on these positive strengths. Thus the almost-positivity result and its exact zero set below hold for this whole family. We keep $s_1=s_2=s_3=1$ for specificity.

The submersion model is the product of $\Sp(3),L_1,L_2,L_3$, all with the inner products induced by $Q$. Right translations on $\Sp(3)$ remain isometries, so its quotient by $H$ is a smooth Riemannian manifold. Proposition~\ref{prop:flatcriterion} proves nonnegative curvature by a quotient construction.

For a skew-Hermitian quaternionic matrix $V$ with diagonal entries $v_1,v_2,v_3\in\Imag\HH$,
\begin{equation}\label{q:eq:diagprojections}
 P_2V=\diag(v_1,v_2,v_3),\qquad
 P_3V=\tfrac13\diag(v_1+v_2+v_3,v_1+v_2+v_3,v_1+v_2+v_3).
\end{equation}
Thus, in addition to the first-block condition, the new conditions are
\begin{equation}\label{q:eq:diagonalconditions}
 [v_j,w_j]=0\ (j=1,2,3),\qquad
 [v_1+v_2+v_3,w_1+w_2+w_3]=0.
\end{equation}
Two imaginary quaternions commute if and only if their real three-dimensional vectors are parallel. This elementary fact is the new restriction unavailable for the abelian diagonal algebra of $\mathfrak{su}(3)$.

\subsection{A complete parametrization of the zero locus}
\subsubsection{The parameter sets}
Use the fixed $F$ and $\NN$ from \eqref{q:eq:Fourier} and \eqref{q:eq:normalizer}. For $0\leq\theta,\beta\leq\pi/2$ put
\begin{equation}\label{q:eq:R}
 R_{\theta,\beta}=
 \begin{pmatrix}
 c_\theta&-s_\theta c_\beta&s_\theta s_\beta\\
 s_\theta&c_\theta c_\beta&-c_\theta s_\beta\\
 0&s_\beta&c_\beta
 \end{pmatrix},\qquad c_u=\cos u,\quad s_u=\sin u.
\end{equation}
This is a real special orthogonal matrix. Let $P_{\theta,\beta}$ be its matrix of entrywise squares:
\begin{equation}\label{q:eq:P}
 P_{\theta,\beta}=
 \begin{pmatrix}
 c_\theta^2&s_\theta^2c_\beta^2&s_\theta^2s_\beta^2\\
 s_\theta^2&c_\theta^2c_\beta^2&c_\theta^2s_\beta^2\\
 0&s_\beta^2&c_\beta^2
 \end{pmatrix}.
\end{equation}
For $\epsilon=(\epsilon_1,1,\epsilon_3)$ with $\epsilon_1,\epsilon_3\in\{\pm1\}$, let
\[
 J_\epsilon=\diag(j_1,1,j_3),\qquad
 j_r=\begin{cases}1,&\epsilon_r=1,\\\jj,&\epsilon_r=-1.\end{cases}
 \qquad
 \Lambda=\begin{pmatrix}1&0\\-1&1\\0&-1\end{pmatrix}.
\]
For $D=\diag(d_1,d_2,d_3)\in K$, define $N_D$ to be the real $3\times3$ matrix with columns the coordinates, in the basis $\ii,\jj,\kk$, of
\begin{equation}\label{q:eq:axes}
 n_r=d_r\ii\bar d_r\in S^2\subset\Imag\HH.
\end{equation}
Define a compact set by the following explicit parametrization:
\begin{equation}\label{q:eq:paramset}
\begin{split}
 \PP=\{\,&D R_{\theta,\beta}J_\epsilon\,nF^*H:
 D\in\Sp(1)^3,\ n\in\NN,\ 0\leq\theta,\beta\leq\pi/2,\\
 &\epsilon=(\pm1,1,\pm1),\quad
 \rank\bigl(N_D P_{\theta,\beta}\diag(\epsilon)\Lambda\bigr)\leq1\,\}.
\end{split}
\end{equation}
The rank requirement means that the three $2\times2$ minors of the displayed real $3\times2$ matrix vanish. Thus \eqref{q:eq:paramset} contains no curvature quantities that still need to be computed. All parameter ranges and matrix entries are explicit. The parametrization is not intended to be one-to-one.

\begin{theorem}[Exact zero-point locus]\label{q:thm:exact}
For the metric \eqref{q:eq:weights} on $G/H$, with $H=\Sp(1)^3\rtimes\Gamma$ and $\Gamma=C_3$ or $S_3$,
\begin{equation}\label{q:eq:exactzero}
 {\ZZ(g_*)=\BB\cup\PP,\qquad
 \BB=\{aH:\Psi(aH)=\tfrac12\}.}
\end{equation}
This is a complete compact parametrization of the zero-point locus. Unlike \eqref{eq:oldall}, it is not asserted to be a single inequality in the last-row norms.
\end{theorem}

\subsubsection{Reduction after the first block condition}
\begin{lemma}\label{q:lem:blockform}
Every plane satisfying the first-block condition has, after right multiplication of its eigenbasis by an element of $\NN$, an eigenbasis of the form
\begin{equation}\label{q:eq:intermediateB}
 B=\diag(T,d_3)R_{0,\beta}J_{(1,1,\epsilon_3)},
 \qquad T\in\Sp(2),\ d_3\in\Sp(1),\quad 0\leq\beta\leq\pi/2.
\end{equation}
The two upper-block eigenvalue forms are independent except when
\begin{equation}\label{q:eq:exception}
 \epsilon_3=1,\qquad \beta=\pi/4.
\end{equation}
In that exceptional case the projection of the entire plane onto $\lfr_1$ has dimension one.
\end{lemma}
\begin{proof}
Lemma~\ref{q:lem:row} permits a column permutation so that the first last-row entry vanishes. If the other two are nonzero, their ratio belongs to $\C^\times$ or to $\C^\times\jj$, respectively, by \eqref{q:eq:axisparallel}. Complex phases in the columns and a left unit quaternion in the last row put that row in the form
$(0,d_3\sin\beta,d_3\cos\beta\,j_3)$. A row with one nonzero entry is covered by an endpoint. Orthogonality then makes the upper two rows a quaternionic orthonormal basis of the complement, giving \eqref{q:eq:intermediateB}.

For $\lambda_1+\lambda_2+\lambda_3=0$, the upper-left block of $B(i\diag\lambda)B^*$ is
\begin{equation}\label{q:eq:uppercompression}
 T\diag\bigl(i\lambda_1,
 i(c_\beta^2\lambda_2+\epsilon_3s_\beta^2\lambda_3)\bigr)T^*.
\end{equation}
The bottom entry is
$d_3i(s_\beta^2\lambda_2+\epsilon_3c_\beta^2\lambda_3)\bar d_3$.
On the trace-zero plane, the coefficient of $\lambda_2$ in the second form of \eqref{q:eq:uppercompression} is $c_\beta^2-\epsilon_3s_\beta^2$. It vanishes precisely in \eqref{q:eq:exception}. There the two blocks together depend only on $\lambda_1$, so the whole $\lfr_1$-projection has rank one.
\end{proof}

\subsubsection{The exceptional case is the entire one-step boundary}
\begin{lemma}\label{q:lem:boundary}
All points represented by the exceptional case \eqref{q:eq:exception} lie in $\BB$, and every point of $\BB$ has such a representation. All corresponding planes satisfy every deeper projection condition, including those for $L_2$ and $L_3$.
\end{lemma}
\begin{proof}
A last row with two complex entries of squared norms $u,1-u$ has, after multiplication by $F^*$, fourth moment
\begin{equation}\label{q:eq:fourthmoment}
 \sum_{j=1}^3 |(bF^*)_j|^4
 =\frac13+\frac23u(1-u).
\end{equation}
This is independent of their complex phases: in the expansion, the nonconstant Fourier terms sum to zero. Permutations and the common $\jj$ in \eqref{q:eq:normalizer} do not change the formula. A common left quaternion preserves norms. In the exceptional case $u=1/2$, so \eqref{q:eq:fourthmoment} equals $1/2$.

Conversely, the construction in the proof of Theorem~\ref{thm:loci} gives, at any point with $\Psi=1/2$, a complex last row with one zero entry. Formula \eqref{q:eq:fourthmoment} forces the other entries to have equal norm. Complex eigenbasis phases put them in the form required by \eqref{q:eq:exception}. Finally, a rank-one $\lfr_1$-projection has rank at most one after any further linear projection, so all additional brackets vanish.
\end{proof}

\subsubsection{The nonexceptional case}
\begin{lemma}[A two-row quaternionic calculation]\label{q:lem:U2}
Suppose $T\in\Sp(2)$ and
\begin{equation}\label{q:eq:rowaxes}
 [T_{r1}i\overline{T_{r1}},T_{r2}i\overline{T_{r2}}]=0
 \quad(r=1,2).
\end{equation}
Then
\begin{equation}\label{q:eq:Tform}
 T=\diag(d_1,d_2)\,U\,\diag(j_1,1),
 \qquad U\in\UU(2),\quad j_1\in\{1,\jj\}.
\end{equation}
The statement includes zero entries of $T$.
\end{lemma}
\begin{proof}
If both entries in the first row are nonzero, \eqref{q:eq:rowaxes} says that their ratio is in $\C^\times$ or $\C^\times\jj$. A left unit quaternion and, if necessary, multiplication of the second column by $-\jj$ make the first row complex. Its quaternionic orthogonal complement is the left-quaternionic span of its complex unit perpendicular row. Multiplying that second row by a left unit quaternion therefore makes the whole matrix complex unitary. When an entry vanishes the same argument applies directly.

This initially gives a factor $\diag(1,\jj)$ on the right. It may be moved to the first column using, for complex $U$,
\[
 U\diag(1,\jj)=\diag(\jj,\jj)\,\bar U\,\diag(-\jj,1).
\]
Absorb the complex minus sign and the left diagonal factor to obtain \eqref{q:eq:Tform}.
\end{proof}

\begin{lemma}\label{q:lem:normalformB}
Outside the exceptional case \eqref{q:eq:exception}, the first two projection conditions hold if and only if an eigenbasis can be written, up to right multiplication by $\NN$, as
\begin{equation}\label{q:eq:Bnormal}
 B=D R_{\theta,\beta}J_\epsilon
\end{equation}
with $D,\theta,\beta,\epsilon$ in the ranges specified in \eqref{q:eq:paramset}, but without imposing its final rank condition.
\end{lemma}
\begin{proof}
The two forms in \eqref{q:eq:uppercompression} are independent. Hence the vanishing of the commutators of the first two diagonal entries is exactly \eqref{q:eq:rowaxes}; the third diagonal entry already lies on a fixed imaginary axis. Apply Lemma~\ref{q:lem:U2}.

Every complex $2\times2$ unitary matrix is a product of two diagonal unitary matrices and a real rotation $R_\theta$, $0\leq\theta\leq\pi/2$: choose the phases of its first row, then use the orthogonality of the second. Absorb the left phases into $d_1,d_2$. The right phases can be moved through the $2,3$ rotation by
\begin{equation}\label{q:eq:phaseidentity}
 \diag(z_1,z_2,1)R_{0,\beta}
 =\diag(1,1,\bar z_2)R_{0,\beta}\diag(z_1,z_2,z_2),
 \qquad |z_1|=|z_2|=1.
\end{equation}
The first-column $\jj$ commutes with that real $2,3$ rotation. Moving complex phases past a $\jj$ merely conjugates them. Absorbing the last-row phase into $d_3$ and the final column phases into $\NN$ gives \eqref{q:eq:Bnormal}. Conversely its two upper-block diagonal axes satisfy \eqref{q:eq:rowaxes}, and its last row satisfies Lemma~\ref{q:lem:row}. Thus both projection conditions hold.
\end{proof}

For $B$ in \eqref{q:eq:Bnormal}, direct multiplication gives the diagonal entries of $B(i\diag\lambda)B^*$ as
\begin{equation}\label{q:eq:diagonalcoeffs}
 n_r\bigl(P_{\theta,\beta}\diag(\epsilon)\lambda\bigr)_r,
 \qquad r=1,2,3.
\end{equation}
Their sum, as a real three-vector, is therefore
$N_D P_{\theta,\beta}\diag(\epsilon)\lambda$.
The final bracket condition in \eqref{q:eq:diagonalconditions} is exactly the rank condition in \eqref{q:eq:paramset}.

\begin{proof}[Proof of Theorem~\ref{q:thm:exact}]
Take any zero plane. Proposition~\ref{prop:flatcriterion} and Corollary~\ref{q:cor:eigenplane} give its normal eigenplane. Lemma~\ref{q:lem:blockform} applies. In its exceptional case Lemma~\ref{q:lem:boundary} places the base point in $\BB$. Otherwise Lemma~\ref{q:lem:normalformB}, \eqref{q:eq:diagonalcoeffs}, and \eqref{q:eq:basefromB} place it in $\PP$.

Conversely a point of $\BB$ has a flat plane by Lemma~\ref{q:lem:boundary}. At a point parametrized in \eqref{q:eq:paramset}, the plane $\Ad_B\afr$ is normal horizontal by \eqref{q:eq:basefromB}; its first two projection brackets vanish by Lemma~\ref{q:lem:normalformB}, or directly from the matrices at exceptional parameters. Equation \eqref{q:eq:diagonalcoeffs} gives the last vanishing bracket. Sufficiency in Proposition~\ref{prop:flatcriterion} proves that the quotient plane is flat.
\end{proof}

\begin{proposition}[The entire zero-plane locus]\label{q:prop:allplanes}
All zero planes are obtained as follows. Use either
\begin{equation}\label{q:eq:twofamilies}
\begin{array}{ll}
 B=D R_{\theta,\beta}J_\epsilon
 &\text{with the rank condition in \eqref{q:eq:paramset}},\\[2pt]
 B=\diag(T,d_3)R_{0,\pi/4}
 &\text{with }T\in\Sp(2),\ d_3\in\Sp(1).
\end{array}
\end{equation}
Choose $n\in\NN$, set $a=BnF^*$, and put
\[
 V=B\diag(i,-i,0)B^*,\qquad W=B\diag(0,i,-i)B^*.
\]
At $aH$ project the plane spanned by
$dR_a\Phi_*^{-1}V$ and $dR_a\Phi_*^{-1}W$.
Every plane so produced has zero curvature, and every zero plane occurs. Changing the representative of $aH$ transports the plane by the quotient identification.
\end{proposition}
\begin{proof}
Both eigenplane families were proved sufficient above. Conversely the reduction of an arbitrary zero plane was a reduction of the whole two-dimensional subspace, not just of one of its vectors. The bases displayed span $\afr$, so the proof of Theorem~\ref{q:thm:exact} and \eqref{eq:projectedplane} give precisely the asserted enumeration.
\end{proof}

\subsubsection{Three explicit checks on the locus}
At $eH$ one has $\Psi=1$, so every two-plane is positive for $g_*$. At $F^*H$, take $B=I$ and the two diagonal generators in Proposition~\ref{q:prop:allplanes}. They satisfy every bracket equation, so $F^*H$ is a zero point with $\Psi=1/3$. This shows that the new zero locus is not merely the boundary $\BB$.

For an explicit regular boundary point, take
\[
 a=R_{0,\pi/4}F^*.
\]
Its last row is $(\sqrt{2/3},-1/\sqrt6,-1/\sqrt6)$, so the squared norms are $(2/3,1/6,1/6)$ and $\Psi=1/2$. For the plane with eigenbasis $B=R_{0,\pi/4}$, the generator $B\diag(0,i,-i)B^*$ has zero $\lfr_1$-projection. Thus it is a flat plane for every choice of positive strengths in \eqref{q:eq:family}, not just for the numerical weights in \eqref{q:eq:weights}.

\subsection{The Hausdorff dimension of the zero locus}
An explicit zero-locus parametrization is not yet a proof of almost positivity: an image of a compact parameter space can have interior. We now bound its dimension, taking account of continuous phase redundancies.

\subsubsection{The last-row strata}
Put $r_j=|a_{3j}|^2$. On the ordered cover $G/K$, the map
\[
 r=(r_1,r_2,r_3):G/K\longrightarrow
 \{r_j\geq0,\ r_1+r_2+r_3=1\}
\]
has rank two whenever all $r_j>0$. Indeed the last-row map $\Sp(3)\to S^{11}$ is a submersion with fiber $\Sp(2)$, and varying the three positive norms on the sphere realizes both independent probability variations. These descend through $K$. The finite quotient does not affect local dimensions, by \eqref{eq:coverZ}.

The circle $\sum r_j^2=1/2$ is regular in the interior of that simplex. Its preimage has dimension $11$. Where that circle touches an edge it has $r=(1/2,1/2,0)$, up to permutation. The corresponding row space has dimension six, so its preimage in $\Sp(3)$ has dimension $6+10=16$; quotienting by $K$, of dimension nine, gives dimension seven. Consequently
\begin{equation}\label{q:eq:boundarydimension}
 \dim_{\mathrm H}\BB=11,\qquad\Vol(\BB)=0.
\end{equation}
The stratum
\begin{equation}\label{q:eq:barycenter}
 \CC=\{r_1=r_2=r_3=1/3\}
\end{equation}
has dimension ten. We use $\CC$ as a containing set in the estimates below; we do not assert that every point of $\CC$ is a zero point of $g_*$.

If $\epsilon_3=-1$ in \eqref{q:eq:paramset}, the two nonzero entries of the last row of $B$ have opposite imaginary axes. After complex eigenphases, their products with the entries of $F^*$ are in orthogonal real two-planes $\C$ and $\C\jj$. Their squared norms therefore add without a cross term, giving
\begin{equation}\label{q:eq:oppositebary}
 |(BnF^*)_{3j}|^2=1/3\quad(j=1,2,3).
\end{equation}
Hence this entire parameter branch lies in $\CC$. The endpoints $\beta=0,\pi/2$ have the same property. For $\epsilon_3=1$, formula \eqref{q:eq:fourthmoment} gives
\begin{equation}\label{q:eq:betaPsi}
 \Psi=\frac13+\frac16\sin^2(2\beta).
\end{equation}
Thus $\beta=\pi/4$ lies in $\BB$. It remains to estimate only
\begin{equation}\label{q:eq:remainingbranch}
 \epsilon_3=1,\qquad
 0<\beta<\pi/2,\qquad\beta\neq\pi/4.
\end{equation}

\subsubsection{The rank condition on the three axes}
\begin{lemma}[A real three-vector estimate]\label{q:lem:axisdim}
Fix a rank-two real $3\times2$ matrix $C_0$. Among triples
$(n_1,n_2,n_3)\in(S^2)^3$, the condition
\begin{equation}\label{q:eq:axisrank}
 \rank\bigl([\,n_1\ n_2\ n_3\,]C_0\bigr)\leq1
\end{equation}
has dimension at most four. More precisely it is covered by countably many smooth images of subsets of at most four-dimensional coordinate domains. If $C_0$ varies smoothly in $q$ parameters while retaining rank two, the bound is $4+q$.
\end{lemma}
\begin{proof}
Choose two independent rows of $C_0$, permute their labels, and multiply on the right by an invertible matrix. Its rows become $(1,0),(0,1),(a,b)$. The columns in \eqref{q:eq:axisrank} are $n_1+a n_3$ and $n_2+b n_3$.

When $n_1,n_2$ are dependent, choose $n_1\in S^2$, its sign for $n_2$, and $n_3\in S^2$, giving at most four parameters. Otherwise their cross product is nonzero. If $(a,b)=(0,0)$, no solution exists. In the other cases the vanishing cross product of the two columns first forces $n_3$ into the plane spanned by $n_1,n_2$: its component normal to that plane would contribute the nonzero vector $(b n_1-a n_2)\times n_3^\perp$. Write $n_3=xn_1+yn_2$. The remaining equations are
\begin{equation}\label{q:eq:ellipse}
 ax+by=-1,\qquad x^2+y^2+2(n_1\cdot n_2)xy=1.
\end{equation}
A line intersects this positive-definite ellipse in at most two points. The four parameters $n_1,n_2$ therefore suffice.

For the stated smooth-image assertion, solve the line equation and then its quadratic intersection with the ellipse. Away from a tangent intersection each of the two roots depends smoothly on the coefficients. At a tangent intersection its double root is the smooth midpoint formula of the quadratic. Work in countably many charts where the denominators used are nonzero. This covers all solutions, including tangent ones, by smooth maps defined on subsets of the coefficient domains. The same construction permits the additional $q$ parameters of $C_0$.
\end{proof}

A direct determinant calculation in \eqref{q:eq:P} gives
\begin{equation}\label{q:eq:detP}
 \det P_{\theta,\beta}=\cos(2\theta)\cos(2\beta).
\end{equation}
On \eqref{q:eq:remainingbranch}, when $\theta\neq\pi/4$, the matrix
$P_{\theta,\beta}\diag(\epsilon)\Lambda$ has rank two. Lemma~\ref{q:lem:axisdim} bounds its allowed axes and two angles together by dimension six.

\subsubsection{Removing phase redundancy}
\begin{lemma}[Four, not six, independent phase parameters]\label{q:lem:phases}
For fixed axes $n_1,n_2,n_3$, fixed angles, fixed $\epsilon$, and a fixed component of $\NN$, the base-point map in \eqref{q:eq:paramset} requires at most four real phase parameters.
\end{lemma}
\begin{proof}
On local charts of the Hopf map $\Sp(1)\to S^2$, write
$d_r=s(n_r)e^{i\alpha_r}$. The three eigenphases in $n=EP$ or $EP(\jj I_3)$ give three more phases. There are two independent redundancies.

First, multiplying all three eigenphases by the same $e^{it}$ multiplies the representative on the right by $e^{it}I_3$ or $e^{-it}I_3$, after moving it through the permutation, possible common $\jj$, and $F^*$. This is a right $K$ action and does not change the point.

Second, replacing all $d_r$ by $d_re^{it}$ leaves their axes fixed. As $R_{\theta,\beta}$ is real,
\[
 e^{it}I_3\,R_{\theta,\beta}J_\epsilon
 =R_{\theta,\beta}J_\epsilon\diag(e^{i\epsilon_1t},e^{it},e^{i\epsilon_3t}).
\]
Compensating by the inverse diagonal factor in $E$ preserves the representative. Use this freedom to set $\alpha_3=0$, then the first freedom to set one eigenphase to zero. These operations are independent even when all signs are positive. The four remaining phases, local Hopf charts, and finitely many discrete components cover the whole parameter map.
\end{proof}

\begin{proof}[Proof of Theorems~\ref{thm:existence} and~\ref{thm:loci} for $d=4$]
Nonnegativity follows from Proposition~\ref{prop:flatcriterion}. The exact locus was established in Theorem~\ref{q:thm:exact}. Its branches with $\epsilon_3=-1$, or $\beta=0,\pi/2$, lie in the ten-dimensional set $\CC$. The branch $\epsilon_3=1$, $\beta=\pi/4$, and the additional exceptional family, lie in $\BB$, of dimension eleven.

Consider \eqref{q:eq:remainingbranch}. If $\theta=\pi/4$, even without imposing the rank equation we have at most six axis parameters, one angle, and four phases: dimension at most eleven. If $\theta\neq\pi/4$, Lemma~\ref{q:lem:axisdim} and \eqref{q:eq:detP} give at most six parameters for the allowed axes and two angles, then Lemma~\ref{q:lem:phases} adds at most four phases: dimension at most ten. Thus the entire zero set is covered by countably many smooth images of subsets of spaces of dimension at most eleven.

The dimension rules of Section~\ref{sec:prelim} now give
$\dim_{\mathrm H}\ZZ(g_*)\le11$ and $\Vol_{g_*}\ZZ(g_*)=0$.

The lower bound eleven follows from the regular part of $\BB\subset\ZZ(g_*)$, for instance above $r=(1/6,1/6,2/3)$. The zero set is closed, either by the explicit compact parametrization or by compactness of the plane bundle and continuity of curvature. A closed set of volume zero has empty interior. Its complement is therefore open and dense, and every plane based there is strictly positive. Finally, the first-block criterion is a necessary part of the new flat equations, so \eqref{eq:oldall} proves $\ZZ(g_*)\subset\{\Psi\leq1/2\}$. This proves \eqref{eq:mainloci} for $d=4$; the topological assertions were established in Section~\ref{sec:topology}.
\end{proof}

\begin{remark}
The parameter count is not merely a count of equations. Lemma~\ref{q:lem:axisdim} controls their actual solution set, and Lemma~\ref{q:lem:phases} exhibits the two genuine continuous redundancies. Ignoring either issue would not prove density.
\end{remark}

\section{Exceptional flags and normal commuting planes}\label{sec:exceptional}\label{e:sec:algebra}
Let $J=\operatorname{Herm}_3(\OO)$ be the Albert algebra, with Jordan product $A\circ B=(AB+BA)/2$ and inner product $\langle A,B\rangle=\tr(A\circ B)$. Its trace-one primitive idempotents form $\OP^2$. Fix the coordinate idempotents $E_1,E_2,E_3$. We use the standard identifications~\cite{Baez,Kollross}
\[
 G=\Aut(J)=F_4,\quad L_9=\Stab_G(E_3)=\Spin(9),
 \quad K=\Stab_G(E_1,E_2,E_3)=\Spin(8).
\]
Let $H=K\rtimes\Gamma$ for $\Gamma=C_3$ or $S_3$, and put $W=G/K=W_8$ and $\pfr=\mathfrak l_9^\perp$. Write $g_*=g_{\OO}$ for the exceptional metric. The vector representation of $L_9$ is its action on $V=\operatorname{Herm}_2(\OO)_0\cong\R^9$ in the upper block. Choose an orthonormal basis $e_1,\ldots,e_9$ with $e_9$ proportional to $E_1-E_2$. For $3\le j\le8$, define $L_j$ to be the connected lift of the subgroup $\SO(j)$ acting on $\Span(e_1,\ldots,e_j)$ and fixing the orthogonal complement. Then $L_8=K$, and these embeddings specify \eqref{eq:exceptionalchain} completely.

For unit strengths, $\Phi_*^{-1}=I+\sum_{j=3}^9P_j$. If $\mathfrak v_j=\mathfrak l_j\ominus\mathfrak l_{j-1}$ for $4\le j\le9$, the right-invariant metric is
\begin{equation}\label{e:eq:metric}
 \overline g_e=Q|_{\pfr}
 +\frac12Q|_{\mathfrak v_9}+\frac13Q|_{\mathfrak v_8}
 +\frac14Q|_{\mathfrak v_7}+\frac15Q|_{\mathfrak v_6}
 +\frac16Q|_{\mathfrak v_5}+\frac17Q|_{\mathfrak v_4}
 +\frac18Q|_{\mathfrak l_3}.
\end{equation}
The summand dimensions are $16,8,7,6,5,4,3,3$. Proposition~\ref{prop:flatcriterion} applies; all strengths can be taken arbitrarily positive without changing the zero-locus equations. We retain the notation
\begin{equation}\label{e:eq:rPsi}
 r_i(aK)=\langle E_3,aE_i\rangle,\qquad \Psi=\sum r_i^2.
\end{equation}

\subsection{Triality data}

We use the octonion norm, conjugation, and alternativity. In particular
\[
 x\bar x=|x|^2,\qquad L_xL_{\bar x}=|x|^2I,
 \qquad [x,y,z]=(xy)z-x(yz)\text{ is alternating}.
\]
Fix a unit imaginary octonion $i$ and the associative subalgebra $\C=\Span_\R\{1,i\}\subset\OO$.

The Lie algebra has the triality model
\begin{equation}\label{e:eq:trialitydecomp}
 \gfr=\frso(8)\oplus\OO\oplus\OO\oplus\OO,
 \qquad \kfr=\frso(8),\qquad \kfr^\perp=\OO^3.
\end{equation}
The three octonion summands are the three inequivalent eight-dimensional $K$-modules. Let $\lambda$ be the order-three triality automorphism of $\frso(8)$. With $u\wedge v=uv^T-vu^T$, the part of the bracket needed for two horizontal vectors is
\begin{align}\label{e:eq:horizontalbracket}
 [(0,u,v,w),(0,x,y,z)]
  =\big(&-4(u\wedge x+\lambda^2(v\wedge y)+\lambda(w\wedge z)),\notag\\
        &\overline{vz-yw},\ \overline{wx-zu},\ \overline{uy-xv}\big).
\end{align}
This is the compact-real-form bracket of Kollross~\cite[Theorem 3.1]{Kollross}. We need only the following special triality identities. For the fixed imaginary unit $i$ and $E=1\wedge i$,
\begin{equation}\label{e:eq:trialityidentity}
 \lambda E=\tfrac12L_i,\qquad \lambda^2E=\tfrac12R_i,
 \qquad E+\lambda E+\lambda^2E=0.
\end{equation}
Indeed $L_i+R_i=-2E$. Moreover $E$ and $\lambda E$ are linearly independent: the first has rank two and the second rank eight. These formulas use the orthogonal pair $1,i$; no formula for arbitrary nonorthogonal decomposable inputs is needed.

The triality action preserves the real trilinear form
\begin{equation}\label{e:eq:trilinear}
 T(u,v,w)=\Real((uv)w).
\end{equation}
The triality action of $K=\Spin(8)$ is transitive on the product of unit spheres in any two of its inequivalent eight-dimensional modules, and the stabilizer of $(1,1)$ is $G_2$ acting simultaneously on all three modules. Also $G_2$ is transitive on the imaginary unit sphere, with stabilizer $\SU(3)$~\cite[Sections 3--4]{Baez}\cite[Sections 2--3]{Kollross}. We record the first fact as
\begin{equation}\label{e:eq:pairtransitive}
 K\text{ acts transitively on }S^7\times S^7.
\end{equation}

\subsection{The rank-one symmetric complement}
The symmetric decomposition for $G/L_9=\OP^2$ is
\begin{equation}\label{e:eq:symmetric}
 \gfr=\lfr_9\oplus\pfr,\qquad
 [\lfr_9,\pfr]\subset\pfr,\quad [\pfr,\pfr]\subset\lfr_9,
 \qquad \pfr\cong\OO^2.
\end{equation}
The $L_9$-representation on $\pfr$ is the real spin representation, transitive on its unit sphere~\cite[Section 4.2]{Baez}. Equivalently, the compact rank-one symmetric space $\OP^2$ has positive sectional curvature. We will use the standard consequence
\begin{equation}\label{e:eq:rankone}
 U,V\in\pfr,\quad [U,V]=0
 \quad\Longrightarrow\quad U,V\text{ are real-linearly dependent}.
\end{equation}
Triality gives the same assertion for the sum of any two of the three octonion modules in \eqref{e:eq:trialitydecomp}. For comparison with the mixed condition \eqref{eq:fatness}, the action of $x\in\OO\cong\mathfrak{spin}(9)\ominus\mathfrak{spin}(8)$ on $\OO\oplus\OO$ is
\[
 T_x=\frac12\begin{pmatrix}0&L_x\\-L_{\bar x}&0\end{pmatrix},
 \qquad T_x^2=-\frac{|x|^2}{4}I.
\]
Thus every nonzero vertical vector acts invertibly on the horizontal spin module. This also verifies \eqref{eq:fatness} directly in the exceptional case.

\subsection{Classification of normal horizontal commuting planes}\label{e:sec:normalform}

In the triality coordinates, put
\begin{equation}\label{e:eq:q0}
 X_0=(0,1,1,1),\qquad Y_0=(0,i,i,i),\qquad
 \qfr_0=\Span_\R\{X_0,Y_0\}\subset\kfr^\perp.
\end{equation}
Equations \eqref{e:eq:horizontalbracket}--\eqref{e:eq:trialityidentity} show that $[X_0,Y_0]=0$.

\begin{lemma}\label{e:lem:normalform}
Every two-dimensional abelian subspace of $\kfr^\perp$ is $\Ad_k\qfr_0$ for some $k\in K$.
\end{lemma}
\begin{proof}
Let independent commuting vectors be $X=(0,u,v,w)$ and $Y=(0,x,y,z)$. None of the three component projections vanishes identically on their plane: otherwise the plane lies in the sum of two modules, contradicting \eqref{e:eq:rankone}. Choose $X$ with all three components nonzero. By \eqref{e:eq:pairtransitive}, arrange $u=a$, $v=b$ with $a,b>0$ real.

The octonionic components of $[X,Y]=0$ give
\[
 bz=yw,\qquad wx=az,\qquad ay=bx.
\]
Thus $y=(b/a)x$, $z=xw/a$, and $wx=xw$. Subtract a real multiple of $X$ from $Y$ to make $x$ imaginary. It remains nonzero, since otherwise all of $Y$ vanishes. Write $x=ri$, after choosing the imaginary unit in this argument, with $r\ne0$. Commuting octonions have parallel imaginary parts, so $w\in\Span\{1,i\}$. The $\frso(8)$ component of the bracket is now, up to the nonzero factor $-4r/a$,
\[
 a^2E+b^2\lambda^2E+|w|^2\lambda E=0.
\]
Equation \eqref{e:eq:trialityidentity} implies $a^2=b^2=|w|^2$. Rescale the two spanning vectors to obtain
\[
 X=(0,1,1,e^{i\phi}),\qquad Y=(0,i,i,ie^{i\phi}).
\]
Replace them by the rotated basis $X_\theta=\cos\theta X+\sin\theta Y$, $Y_\theta=-\sin\theta X+\cos\theta Y$. The three components of $X_\theta$ are $e^{i\theta},e^{i\theta},e^{i(\phi+\theta)}$, so their invariant trilinear product has real part $\cos(\phi+3\theta)$. Choose $\theta$ to make this $1$. Normalize its first two unit components to $1,1$ by $K$; invariance of \eqref{e:eq:trilinear} forces its last component to be $1$ as well. The off-module bracket equations then force the three components of the transformed $Y_\theta$ to coincide. They are unit imaginary octonions, by componentwise orthogonality. The common stabilizer $G_2$ carries this unit to the fixed $i$. This gives \eqref{e:eq:q0}.
\end{proof}

\subsubsection{The embedded complex group and a Fourier plane}
For a trace-zero skew-Hermitian complex $3\times3$ matrix $Z$, the operation
\[
 D_Z(A)=ZA-AZ\qquad(A\in J)
\]
is an Albert-algebra derivation. This is the standard decomposition
$\Der(J)=\Der(\OO)\oplus\mathfrak{sa}_3(\OO)$; the second summand acts by matrix commutators~\cite[Section 4.2]{Baez}. Restricted to $Z\in\mathfrak{su}(3)$ these derivations form a Lie subalgebra, since all scalar coefficients lie in one associative complex subalgebra. Its connected subgroup is an embedded $\SU(3)$, denoted $U_c$. On $J_c=\operatorname{Herm}_3(\C)$ it acts by ordinary unitary conjugation.

Let
\begin{equation}\label{e:eq:torus}
 \afr=\{D_{i\diag(\lambda_1,\lambda_2,\lambda_3)}:
       \lambda_i\in\R,\ \lambda_1+\lambda_2+\lambda_3=0\},
 \qquad T=\exp\afr.
\end{equation}
With $\omega=e^{2\pi i/3}$, take the Fourier matrix
\[
 F=\frac1{\sqrt3}\begin{pmatrix}1&1&1\\1&\omega&\omega^2\\1&\omega^2&\omega\end{pmatrix}.
\]
Multiplication by one unit scalar makes its determinant $1$ without changing its conjugation of diagonal matrices. Denote the resulting element of $U_c$ by $f$. Direct multiplication shows that
\begin{equation}\label{e:eq:fourier}
 \Ad_f\afr=\qfr_0.
\end{equation}
For example the corresponding zero-diagonal complex matrices spanning the plane are
\[
 \begin{pmatrix}0&1&-1\\-1&0&1\\1&-1&0\end{pmatrix},\qquad
 \begin{pmatrix}0&i&i\\i&0&i\\i&i&0\end{pmatrix},
\]
and $F$ simultaneously diagonalizes them. To relate the two models explicitly, one can identify the horizontal coordinate $(0,u,v,w)$ with $-D_{M(u,v,w)}$, where
\[
 M(u,v,w)=\begin{pmatrix}
 0&u&-\bar w\\-\bar u&0&v\\w&-\bar v&0
 \end{pmatrix}.
\]
The minus sign makes the cyclic cross brackets agree with \eqref{e:eq:horizontalbracket}. Thus the two displayed complex matrices give exactly the plane \eqref{e:eq:q0}; an overall common scale would not change it.

Combining Lemma~\ref{e:lem:normalform} and \eqref{e:eq:fourier}, every normal horizontal commuting plane at $aH$ has the form
\begin{equation}\label{e:eq:qakf}
 \qfr=\Ad_{akf}\afr\qquad(k\in K).
\end{equation}
This is an enumeration of all such planes, not merely a family of examples.

\subsection{The torus normalizer and its singular points}\label{e:sec:torus}

\subsubsection{Weights and the ten-dimensional normalizer}
For an off-diagonal entry in position $(i,j)$, the action of \eqref{e:eq:torus} is
\begin{equation}\label{e:eq:weightsJ}
 z\in\C:\quad \dot z=i(\lambda_i-\lambda_j)z,
 \qquad
 z\in\C^\perp:\quad \dot z=i(\lambda_i+\lambda_j)z=-i\lambda_kz,
 \quad\{i,j,k\}=\{1,2,3\}.
\end{equation}
The second formula uses $zi=-iz$ for $z\perp\C$. The three root-weight lines $\lambda_i-\lambda_j$ and the three fundamental-weight lines $\lambda_k$ are six distinct real lines on the trace-zero plane. The diagonal part is the fixed subspace of $T$ in $J$.

Let $S=\Aut(\OO)_i\cong\SU(3)$ act on all matrix entries. It fixes $J_c$ pointwise, centralizes $U_c$, and lies in $K$. Moreover
\begin{equation}\label{e:eq:normalizerdim}
 \calN=N_G(\afr),\qquad
 \Lie(\calN)=\mathfrak c_\gfr(\afr)=\afr\oplus\sfr,
 \qquad \dim\calN=10,\quad\dim S=8.
\end{equation}
Here are details avoiding a root-count assumption. A derivation commuting with $T$ sends each $E_i$ into the diagonal fixed subspace; differentiating $E_i\circ E_i=E_i$ then forces that derivative to be zero. Thus it belongs to $\kfr$. In each triality module it preserves $\C$, by the two different weights in \eqref{e:eq:weightsJ}, and rotates it through some infinitesimal angle $\alpha_i$. Invariance of $T(u,v,w)$, evaluated at $(i,1,1)$, gives $\alpha_1+\alpha_2+\alpha_3=0$. Subtract the unique element of $\afr$ with these three rotations. The remaining triality transformation fixes both $1$ and $i$ in every module, so belongs to $\sfr$. Conversely $\afr\oplus\sfr$ centralizes $\afr$. Finally the Lie algebra of a torus normalizer is its centralizer: conjugation on the integral lattice of the torus is discrete.

Every $n\in\calN$ preserves $J_c$. Indeed, it normalizes $S$, the derived connected subgroup of the centralizer in \eqref{e:eq:normalizerdim}, and $J^S=J_c$. It also preserves $J^T=\Span\{E_1,E_2,E_3\}$ and permutes these primitive idempotents. On $J_c$ it therefore acts by monomial unitary conjugation, possibly composed with complex conjugation. This last elementary description follows by examining the three complex off-diagonal matrix-unit spaces: their phases must multiply consistently under the Jordan product, and their choices of conjugation agree. In particular,
\begin{equation}\label{e:eq:unbiased}
 n f^{-1}E_j\in J_c\quad\text{and its three diagonal entries are }1/3
 \qquad(n\in\calN).
\end{equation}

\begin{lemma}[Two parameters for base points]\label{e:lem:basefiber}
Fix a plane $\qfr=\Ad_b\afr$. The set of points $aH$ at which $\qfr\subset(\Ad_a\kfr)^\perp$ has dimension at most two. Explicitly it is
\begin{equation}\label{e:eq:basefiber}
 \{b n f^{-1}H:n\in\calN\}.
\end{equation}
\end{lemma}
\begin{proof}
Equation \eqref{e:eq:qakf} gives \eqref{e:eq:basefiber}, allowing inversion of $n$. Conversely every displayed point is horizontal for the plane, by \eqref{e:eq:fourier}. Since $S\subset K$ centralizes $f$, right multiplication of $n$ by $S$ leaves the point unchanged. The parametrization therefore factors through $\calN/S$, of dimension $10-8=2$.
\end{proof}

\subsubsection{All torus-singular points of the Cayley plane}
Regard $\OP^2$ as the trace-one idempotents of $J$. Let $\calD$ be the set of points fixed by some nonzero element of $\afr$.

\begin{lemma}\label{e:lem:singular}
The set $\calD$ is the union of six spheres: for each pair $\{i,j\}$, one sphere $S^2$ has only its $(i,j)$ off-diagonal entry in $\C$ possibly nonzero, and one sphere $S^6$ has only that entry in $\C^\perp$ possibly nonzero. The remaining diagonal coordinate is zero. Their coordinate endpoints include all $E_i$.
\end{lemma}
\begin{proof}
A nonzero linear functional evaluation on the two-dimensional torus algebra can vanish on at most one of the six distinct weight lines in \eqref{e:eq:weightsJ}. Thus a point fixed by a nonzero torus vector has nonzero entries only in the diagonal and one indicated weight space. For a non-diagonal such idempotent the third diagonal coordinate must be zero. Its remaining block has the form
\[
 \begin{pmatrix}t&z\\\bar z&1-t\end{pmatrix},\qquad |z|^2=t(1-t).
\]
Equivalently $(2t-1)^2+4|z|^2=1$. If $z\in\C$ this is $S^2$; if $z\in\C^\perp$ it is $S^6$. Pure diagonal trace-one idempotents are exactly the $E_i$. The same equations prove the converse.
\end{proof}

\subsubsection{The diagonal map and orbit fibers}
The map $r=(r_1,r_2,r_3)$ in \eqref{e:eq:rPsi} classifies the left $L_9$-orbits on $W$. To see this, pass by inversion to the $K$-action on $\OP^2$. In the chart where the third diagonal entry is positive, a primitive projection has the form
\begin{equation}\label{e:eq:chart}
 x=\frac1{1+|p|^2+|q|^2}
 \begin{pmatrix}|p|^2&p\bar q&p\\qfr\bar p&|q|^2&q\\\bar p&\bar q&1\end{pmatrix},
 \qquad p,q\in\OO.
\end{equation}
The two octonions generate an associative subalgebra, so the formula is unambiguous and is verified directly. Conversely it follows from the rank-one cofactor identities. The $K$-action on $(p,q)$ is the sum of two inequivalent triality modules. Hence \eqref{e:eq:pairtransitive} makes fixed norms a single orbit; the one-zero cases use transitivity on one sphere. Permuting the chart covers the boundary cases.

Varying the two positive norms in \eqref{e:eq:chart} also proves that the diagonal map has rank two in the simplex interior. Since $G\to\OP^2$, $a\mapsto a^{-1}E_3$, is a submersion, the same holds for $r:W\to\R^2$. In particular,
\begin{equation}\label{e:eq:barydimension}
 \mathcal C=\{r_1=r_2=r_3=1/3\}\subset W
 \text{ is a smooth submanifold of dimension }22.
\end{equation}
Its finite-quotient image has the same dimension.

\section{Almost positivity on the exceptional quotients}\label{sec:exceptionalpositive}\subsection{The one-step zero locus}\label{e:sec:first}

Let $\qfr=\Ad_b\afr$ be a commuting plane. The symmetric decomposition and \eqref{e:eq:rankone} give
\begin{equation}\label{e:eq:firstequiv}
 [P_9V,P_9W]=0
 \quad\Longleftrightarrow\quad \dim(\qfr\cap\lfr_9)\ge1
 \quad\Longleftrightarrow\quad b^{-1}E_3\in\calD.
\end{equation}
For the first equivalence, a vanishing $\lfr_9$ bracket and $[V,W]=0$ force the $\pfr$ components to commute, hence to be dependent. Conversely choose $0\ne V\in\qfr\cap\lfr_9$ and separate the $\lfr_9$ and $\pfr$ parts of $[V,W]=0$. The last equivalence says precisely that a nonzero torus vector fixes $b^{-1}E_3$.

\subsubsection{The octonionic singular branch is already small}
If $x=b^{-1}E_3$ belongs to one of the $S^6$ branches of Lemma~\ref{e:lem:singular}, its $J_c$-component is diagonal. For any associated base point $aH=b n f^{-1}H$, \eqref{e:eq:unbiased} gives
\[
 r_i(aH)=\ip{x}{n f^{-1}E_i}=\tfrac13\tr x=\tfrac13.
\]
Thus every zero point arising through these branches belongs to the fixed $22$-dimensional set $\mathcal C$. This includes their coordinate endpoints.

\subsubsection{The complex branch has one rotation parameter}
The $S^2$ branches are complex projective lines. Torus phases and coordinate permutations, all normalizing $\afr$, allow their reduction to
\[
 x_\beta=b_\beta^{-1}E_3,\qquad
 b_\beta=\begin{pmatrix}1&0&0\\0&\cos\beta&-\sin\beta\\0&\sin\beta&\cos\beta\end{pmatrix}\in U_c,
 \qquad 0\le\beta\le\frac\pi2.
\]
Consequently every plane of this branch is
\begin{equation}\label{e:eq:regularq}
 \qfr=\Ad_{l b_\beta}\afr,\qquad l\in L_9.
\end{equation}
On the complex subgroup, $P_9$ is upper-$2\times2$ block projection. Hence
\begin{equation}\label{e:eq:Dbeta}
 P_9\Ad_{b_\beta}(i\diag\lambda)
 =i\diag\big(\lambda_1,c^2\lambda_2+s^2\lambda_3,
                       s^2\lambda_2+c^2\lambda_3\big),
 \quad c=\cos\beta,\ s=\sin\beta.
\end{equation}
The determinant of this map on the trace-zero plane is $\cos(2\beta)$. Therefore
\begin{equation}\label{e:eq:rankdrop}
 P_9\Ad_{b_\beta}\afr=\afr\quad(\beta\ne\pi/4),\qquad
 \dim P_9\Ad_{b_{\pi/4}}\afr=1.
\end{equation}

For $aH=l b_\beta n f^{-1}H$, the monomial/conjugate description of $n|_{J_c}$ gives, after permuting labels,
\begin{equation}\label{e:eq:rphase}
 r_j=\frac{1+\sin(2\beta)\cos(\theta+2\pi(j-1)/3)}3
 \qquad(j=1,2,3)
\end{equation}
for some phase $\theta$. Every phase occurs. Squaring and summing yields
\begin{equation}\label{e:eq:Psimoment}
 \mathord{\Psi=\frac13+\frac16\sin^2(2\beta).}
\end{equation}
Equations \eqref{e:eq:firstequiv}, \eqref{e:eq:barydimension} and \eqref{e:eq:Psimoment} show that the one-step zero set is contained in $\{\Psi\le1/2\}$.

Conversely, the region $r_i\ge0$, $\sum r_i=1$, $\sum r_i^2\le1/2$ is the closed disk inscribed in the probability simplex. Formula \eqref{e:eq:rphase}, with $0\le\sin(2\beta)\le1$ and arbitrary $\theta$, parametrizes this entire disk. Each such model point has a flat plane by \eqref{e:eq:regularq} and Proposition~\ref{prop:flatcriterion}; left $L_9$ invariance of the one-step metric and the orbit statement after \eqref{e:eq:chart} give the same for every point with those diagonal coordinates. This proves \eqref{eq:oldall} for $d=8$. The region $\Psi<1/2$ is nonempty and open.

\subsection{A hypersurface of flats survives every further deformation inside \texorpdfstring{$L_9$}{L9}}\label{e:subsec:boundary}
At $\beta=\pi/4$, \eqref{e:eq:rankdrop} says that every projection of $\qfr$ into a subalgebra of $\lfr_9$ is at most one-dimensional. Thus \emph{all} additional bracket tests vanish. Equation \eqref{e:eq:rphase} now parametrizes the circle $\sum r_i^2=1/2$, and the $L_9$-orbit classification shows
\begin{equation}\label{e:eq:boundaryincluded}
 \mathcal B=\{\Psi=1/2\}\subset\ZZ(g_*).
\end{equation}
This argument does not assert that left $L_9$ remains isometric after all deformations; rather, the rank-one projection property persists after \emph{every} conjugation by $l\in L_9$, so the flat-plane equations themselves hold at all these points.

In the simplex interior the circle is regular and $r$ is a submersion, so $\mathcal B$ has a $23$-dimensional smooth stratum. At an edge the circle touches only $(1/2,1/2,0)$ and its permutations. In $\OP^2$ this fixed-diagonal set is an $S^7$ inside $\OP^1=S^8$. Its inverse image under $G\to\OP^2$ has dimension $7+36$; dividing by $K$, of dimension $28$, gives dimension $15$. Hence
\begin{equation}\label{e:eq:boundarydim}
 \dim_{\mathrm H}\mathcal B=23.
\end{equation}
The remaining task is an upper bound for all the regular complex-branch flats.

\subsection{Commuting compressions of skew-symmetric operators}\label{e:sec:flags}

We give the linear-algebra argument separately because it is the step that upgrades quasipositivity to almost positivity.

For an ordered orthonormal basis of $\R^N$, let $A_j,B_j$ be the compressions of skew-symmetric operators $A,B$ to the span of its first $j$ vectors. We call the basis \emph{admissible} if $[A_j,B_j]=0$ for every $j$. Conditions for $j\le2$ are automatic.

The semialgebraic dimension rules in Section~\ref{sec:prelim} apply to these compression equations, including their singular fibers.

\begin{lemma}[Flag compression estimate]\label{e:lem:flagbound}
Suppose $A,B\in\frso(N)$ satisfy
\[
 \rank A=2,\qquad \rank B=r\ge2,\qquad AB=0.
\]
Then the semialgebraic set of admissible oriented orthonormal bases has dimension at most
\begin{equation}\label{e:eq:flagbound}
 d(N)=1+\frac{(N-1)(N-2)}2.
\end{equation}
In particular the bound is $29$ for $N=9,r=6$.
\end{lemma}
\begin{proof}
Since the operators are skew, $AB=0$ also implies $BA=0$ and $\im A\perp\im B$. Choose the last unit basis vector $e$, and compress to $e^\perp$. If $P$ denotes orthogonal projection, direct multiplication gives
\[
 [PAP,PBP]|_{e^\perp}
 =(Ae)(Be)^T-(Be)(Ae)^T.
\]
The two vectors $Ae,Be$ are orthogonal, so the compressed operators commute exactly when $Ae=0$ or $Be=0$. In that case their product is again zero. The compression of the simple two-form $A$ has rank two unless $e\in\im A$, when it vanishes. The rank of the compression of $B$ is $r-2$ if $e\in\im B$, and $r$ otherwise. The latter assertion follows by splitting off $\ker B$: restriction of a nondegenerate skew form to a hyperplane loses rank two, while a nonzero component of $e$ in its kernel leaves the rank unchanged.

Let $D(N,r)$ be a uniform dimension bound for the admissible bases. If $r=0$, every basis is allowed and $D(N,0)=N(N-1)/2$. For $r\ge2$, the choices of $e$ fall into the following covering of cases:
\begin{center}
\begin{tabular}{p{0.33\textwidth}c p{0.35\textwidth}}
\toprule
Choice of $e$ & Dimension of choices & Remaining operators\\
\midrule
$e\in\im A$ & $1$ & $A$ vanishes; no further restriction\\
$e\in\ker B\setminus\im A$ & $N-r-1$ & ranks $2,r$ on $\R^{N-1}$\\
$e\in\ker A\setminus\im B$ & $N-3$ & ranks $2,r$ on $\R^{N-1}$\\
$e\in\im B$ & $r-1$ & ranks $2,r-2$ on $\R^{N-1}$\\
\bottomrule
\end{tabular}
\end{center}
The middle two cases are absent unless $N-r\ge3$. Some cases overlap, which only makes the resulting upper bound more generous. Their dimensions, together with the fiber rule, give
\begin{align}\label{e:eq:recurrence}
 D(N,r)\le\max\big\{&1+\tfrac{(N-1)(N-2)}2,\notag\\
 &(N-r-1)+D(N-1,r),\quad (N-3)+D(N-1,r),\notag\\
 &(r-1)+D(N-1,r-2)\big\},
\end{align}
omitting the absent cases.

Induct on $N$. The first term is $d(N)$. The middle terms are at most $d(N)-1$, using $r\ge2$ and $d(N)-d(N-1)=N-2$. If $r=2$, the last term is $1+D(N-1,0)=d(N)$. If $r\ge4$, the last term is at most $d(N)$ because $r\le N-2$. This proves \eqref{e:eq:flagbound}, uniformly for every pair in every recursive fiber.
\end{proof}

\subsubsection{Application to the vector representation of \texorpdfstring{$\Spin(9)$}{Spin(9)}}
The vector representation on $V=\operatorname{Herm}_2(\OO)_0$ identifies
\[
 V=\R\oplus\C\oplus\C^\perp.
\]
The torus element in \eqref{e:eq:torus} acts trivially on the first real line, by complex rotation of frequency $\lambda_1-\lambda_2$ on $\C$, and by complex rotation of frequency $\lambda_1+\lambda_2$ on $\C^\perp$. These two frequencies are independent. Therefore its plane of skew operators is
\begin{equation}\label{e:eq:ABnormal}
 \afr=\Span\{A,B\},\qquad
 \rank A=2,\quad \rank B=6,\quad AB=0,
\end{equation}
where $A$ is a unit rotation on a two-plane, $B$ is a complex structure on an orthogonal six-plane, and there is one common kernel line.

The metric $Q$ restricted to $\frso(9)$ is a positive multiple of its trace inner product. Thus each $P_j$ there is precisely upper-left $j\times j$ compression in the chosen vector basis. Define
\begin{equation}\label{e:eq:Eflag}
 \calE=\{l\in L_9:
 [P_j\Ad_l A,P_j\Ad_l B]=0\text{ for }j=3,\ldots,8\}.
\end{equation}
The vector covering $\Spin(9)\to\SO(9)$ has finite fibers. Choosing $l$ in \eqref{e:eq:Eflag} is equivalent to choosing an admissible basis for the fixed pair $A,B$ in Lemma~\ref{e:lem:flagbound}. Consequently
\begin{equation}\label{e:eq:Edim}
 \mathord{\dim\calE\le29.}
\end{equation}
Equivalently, among the $26$-dimensional $\SO(9)$ orbit of the plane \eqref{e:eq:ABnormal}, the allowed planes have dimension at most $19$: its stabilizer has dimension
$\dim\SO(2)+\dim\UU(3)=10$. The estimate \eqref{e:eq:Edim}, rather than this reformulation, will be used below.

\subsection{Dimension of the complete zero locus}\label{e:sec:count}

Consider the regular complex branch \eqref{e:eq:regularq} with $\beta\ne\pi/4$. Since $P_9\Ad_{b_\beta}\afr=\afr$, all additional flat-plane conditions are exactly $l\in\calE$. The candidate planes therefore form the image of
\begin{equation}\label{e:eq:planemap}
 \calE\times\big([0,\pi/2]\setminus\{\pi/4\}\big)
 \longrightarrow\Gr_2(\gfr),\qquad
 (l,\beta)\longmapsto\Ad_{l b_\beta}\afr.
\end{equation}

\begin{lemma}\label{e:lem:stabilizer9}
For every $\beta$, the stabilizer in $L_9$ of $\Ad_{b_\beta}\afr$ has dimension at least nine. For $\beta\ne\pi/4$ its right action preserves the $l$-parameter set $\calE$ in \eqref{e:eq:planemap}.
\end{lemma}
\begin{proof}
The group $S\cong\SU(3)$ of dimension eight centralizes $U_c$, hence centralizes $b_\beta$ and $\afr$. It lies in $L_9$. The circle in $T$ given by $\lambda_2=\lambda_3$ fixes $x_\beta$; its conjugate by $b_\beta$ consequently lies in $L_9$ and centralizes the plane $\Ad_{b_\beta}\afr$. Its Lie algebra is independent of $\sfr$, so together these groups have dimension nine. If $h$ stabilizes this plane and $\beta\ne\pi/4$, equivariance of $P_9$ under $L_9$ shows that $\Ad_h\afr=\afr$. Hence $l$ and $lh$ satisfy exactly the same projection-commutation conditions in \eqref{e:eq:Eflag}.
\end{proof}

Every fiber of \eqref{e:eq:planemap} has dimension at least nine, even when $\beta$ is held fixed. Equations \eqref{e:eq:Edim} and \eqref{e:eq:planemap} therefore give
\begin{equation}\label{e:eq:qdimension}
 \dim\{\text{regular complex-branch candidate planes}\}
 \le29+1-9=21.
\end{equation}
Lemma~\ref{e:lem:basefiber} allows at most two further parameters for their base points. Thus
\begin{equation}\label{e:eq:regularbound}
 \mathord{\dim\{\text{regular complex-branch zero points}\}\le21+2=23.}
\end{equation}

The proof has not omitted a rank-deficient branch. Lemma~\ref{e:lem:singular} enumerates every possibility after the first deformation. Its $S^6$ branches lie in $\mathcal C$, of dimension $22$. The complex branch at $\beta=\pi/4$ lies in $\mathcal B$, of dimension $23$. The remaining branch is covered by \eqref{e:eq:regularbound}. These alternatives exhaust every zero plane by Proposition~\ref{prop:flatcriterion} and Lemma~\ref{e:lem:normalform}. We conclude
\begin{equation}\label{e:eq:upper23}
 \dim_{\mathrm H}\ZZ(g_*)\le23.
\end{equation}
Together with \eqref{e:eq:boundaryincluded}--\eqref{e:eq:boundarydim}, this is equality.

Finally $\ZZ(g_*)$ is closed: a convergent sequence of zero points has a convergent subsequence of zero planes in the compact Grassmann bundle. Since its Hausdorff dimension is less than $24$, it has empty interior and zero Riemannian volume. Its complement is therefore open and dense. Nonnegativity was proved in Proposition~\ref{prop:flatcriterion}; outside the zero locus every sectional curvature is strictly positive. This proves the exceptional curvature and zero-locus assertions in Theorems~\ref{thm:existence} and~\ref{thm:loci}.

\subsection{Exact parametrization of every zero plane}\label{e:sec:exact}

Although the dimension proof does not need a scalar elimination of the zero equations, the construction supplies an exact test. Set
\[
 D_1=D_{i\diag(1,-1,0)},\qquad D_2=D_{i\diag(0,1,-1)}.
\]
These are fixed derivations of $J$, not octonionic matrix-group elements.

\begin{proposition}\label{e:prop:exactlocus}
The complete zero-point locus is
\begin{equation}\label{e:eq:exactparam}
 \mathord{\ZZ(g_*)=
 \left\{b f^{-1}H:
 [P_j\Ad_bD_1,P_j\Ad_bD_2]=0
 \text{ for every }j=3,\ldots,9\right\}.}
\end{equation}
For every $b$ in this parameter set, put $a=b f^{-1}$. The corresponding zero plane is the quotient projection of
\begin{equation}\label{e:eq:allzeroplanes}
 \Span\left\{dR_a(\Phi^{-1}\Ad_bD_1),
                 dR_a(\Phi^{-1}\Ad_bD_2)\right\}.
\end{equation}
Every zero plane occurs in this way.
\end{proposition}
\begin{proof}
For any $b$, \eqref{e:eq:fourier} makes $\Ad_b\afr$ normal-horizontal at $b f^{-1}H$; its generators commute. Thus \eqref{e:eq:exactparam} is sufficient by Proposition~\ref{prop:flatcriterion}. Conversely a flat plane at $aH$ has its normal form $\Ad_{akf}\afr$ by \eqref{e:eq:qakf}. Taking $b=akf$ gives $b f^{-1}H=aH$, the required equations, and precisely its spanning plane. This proves necessity and completeness.
\end{proof}

Formula \eqref{e:eq:exactparam} is a compact real-algebraic parametrization (using $G\subset\SO(J)$), not a claim that every point of $\{\Psi<1/2\}$ survives. Indeed such an open zero region is excluded by \eqref{e:eq:upper23}. The distinguished subset $\{\Psi=1/2\}$ consists of flats that no further nested deformation within $L_9$ can remove.

\paragraph{Two explicit checks.}
At $eH$, $(r_1,r_2,r_3)=(0,0,1)$, so there are no zero planes. For a boundary point take
\[
 a=b_{\pi/4}f^{-1}.
\]
Its diagonal coordinates, after a permutation, are $(2/3,1/6,1/6)$, giving $\Psi=1/2$. The plane $\Ad_{b_{\pi/4}}\afr$ has one-dimensional $\lfr_9$ projection, and therefore satisfies every equation in \eqref{e:eq:exactparam}. This checks explicitly that the new metric is not everywhere positive.

\section{Rigidity of the complex construction and conformal changes}\label{sec:rigidity}
\subsection{A restriction on nested deformations}
\begin{theorem}\label{thm:complexrestriction}
Every metric on $\SU(3)/(T^2\rtimes\Gamma)$, $\Gamma=C_3$ or $S_3$, obtained by a finite nested connected-subgroup construction \eqref{eq:generalPhi} has a nonempty open set of zero-curvature points. In particular none is almost positively curved.
\end{theorem}

\subsubsection{The block subgroup}
For $\mathfrak l_1=\mathfrak{su}(2)\oplus\mathfrak u(1)$,
\[
 [V_{\lfr_1},W_{\lfr_1}]
 =[V_{\mathfrak{su}(2)},W_{\mathfrak{su}(2)}].
\]
A commuting pair in $\mathfrak{su}(2)$ is real-linearly dependent. Every further projection of such a pair is still dependent. Every pair in an abelian Lie algebra commutes automatically. Thus deforming further along the diagonal torus, a circle, or the upper-block $\SU(2)$ does not eliminate any of the zero planes already satisfying the block condition. In particular the region $\{\Psi<1/2\}$ persists for the entire nested block construction, with arbitrary positive deformation strengths.

\subsubsection{The irreducible real subgroup}
There is one other nonabelian proper connected subgroup type to consider, the irreducible real $\SO(3)\subset\SU(3)$. Its Lie algebra decomposition is
\[
 \mathfrak{su}(3)=\mathfrak{so}(3)\oplus i\operatorname{Sym}_0(3,\R).
\]
Take $w\in S^2$ close to $w_0=(1,1,1)/\sqrt3$, let $A(w)$ be the real skew matrix $v\mapsto w\times v$, and let $S$ be a real symmetric trace-zero matrix with $Sw=0$. Such $S$ form a smooth two-dimensional vector bundle near $w_0$. Put
\begin{equation}\label{q:eq:SO3pair}
 W(w)=i(I-3ww^T),\qquad Y(w,S)=A(w)+iS.
\end{equation}
These matrices commute: both $A(w)$ and $S$ preserve $w^\perp$, on which $I-3ww^T$ is scalar. Also $W(w)_{\mathfrak{so}(3)}=0$, so the $\mathfrak{so}(3)$-projection bracket vanishes.

At $a=I,w=w_0,S=0$, the matrices in \eqref{q:eq:SO3pair} have zero diagonal. For $a$ near $I$, solve the four real equations
\begin{equation}\label{q:eq:IFTeq}
 \diag(a^*W(w)a)=0,\qquad\diag(a^*Y(w,S)a)=0.
\end{equation}
Each equation has only two independent imaginary diagonal components. The derivative of the first with respect to $w\in T_{w_0}S^2$ is an isomorphism to the trace-zero diagonal plane: it is $-6i\diag(w_0\odot\dot w)$. At $w_0$ the $S$-space is explicitly
\begin{equation}\label{q:eq:Sslice}
 S(u,v)=\begin{pmatrix}u&-u-v&v\\-u-v&v&u\\v&u&-u-v\end{pmatrix}.
\end{equation}
The first two diagonal entries are $(u,v)$, so the derivative of the second equation with respect to $S$ is also an isomorphism. Its $w$-derivative at $S=0$ has zero diagonal. Thus the full four-variable Jacobian is invertible.

The implicit function theorem supplies $w(a),S(a)$ solving \eqref{q:eq:IFTeq} for every $a$ in a neighborhood of $I$. The two matrices remain independent because $Y$ has a nonzero real skew part and $W$ is purely imaginary symmetric. They give normal horizontal commuting planes whose $\mathfrak{so}(3)$ projections commute. Proposition~\ref{prop:flatcriterion} therefore gives an open zero-point locus for every positive strength of the $\SO(3)$ deformation. All its further proper connected subgroups are abelian, so nested deformations below it do not remove these flats.

\subsubsection{Exhaustion of the subgroup alternatives}
\begin{proof}[Proof of Theorem~\ref{thm:complexrestriction}]
Repeated full-$\SU(3)$ steps only rescale the initial metric and change positive coefficients in \eqref{eq:generalPhi}; they add no condition beyond $[V,W]=0$. Consider the first proper connected subgroup with nonabelian Lie algebra, if one exists.

A reducible subgroup of $\SU(3)$ preserves a complex line and lies in a conjugate of $S(\UU(2)\UU(1))$. Its nonabelian Lie algebra is $\mathfrak{su}(2)$ or $\mathfrak u(2)$: projection to $\mathfrak{su}(2)$ must be onto, since a proper Lie subalgebra of the cross-product algebra $\mathfrak{su}(2)\cong\R^3$ has dimension at most one. A graph into the central line has zero graph homomorphism on the perfect algebra $\mathfrak{su}(2)$. This yields exactly the two alternatives. The block argument above therefore applies to this case and to every nested subgroup below it.

For completeness, an irreducible proper connected subgroup has zero-dimensional Lie algebra center by complex Schur's lemma and the trace-zero condition. Its compact semisimple Lie algebra is simple: an irreducible complex representation of a product is a tensor product, and the prime dimension three permits only one nontrivial factor; faithfulness rules out the others. A simple compact Lie algebra of rank at most two and dimension at most eight is $\mathfrak{su}(2)$ or $\mathfrak{su}(3)$. This is the standard classification of compact simple Lie algebras of rank at most two~\cite{Knapp}. The $\mathfrak{su}(3)$ case is not proper. The irreducible three-dimensional representation of $\mathfrak{su}(2)$ is its real adjoint representation~\cite{Knapp}, so this subgroup is conjugate to $\SO(3)$. The preceding implicit-function argument applies.

If all proper subgroups in the chain are abelian, every additional projection bracket is automatic; the normal horizontal planes from Lemma~\ref{q:lem:normalclassification} remain flat at every point. Finally, conjugating the entire chain simply moves its zero locus by a left translation of $G/H$: under this translation the equation $V\perp\Ad_a\kfr$ and all conjugated projection equations transform together. Each case has a nonempty open set of zero points. This proves the theorem.
\end{proof}

\subsection{Persistence of flat planes under conformal changes}
The next restriction is independent of the subgroup classification above.
\begin{proposition}\label{prop:conformal}
Let $(M,g)$ be a compact normal biquotient of a product of compact groups with bi-invariant metrics. If $\widehat g=e^{2f}g$ has nonnegative sectional curvature, every $g$-zero plane is also a $\widehat g$-zero plane. Hence
\[
 \ZZ(g)\subseteq\ZZ(\widehat g).
\]
In particular no such conformal change turns a one-step metric $g_t$ on $X_d$ or $Y_d$ into an almost-positive metric.
\end{proposition}
\begin{proof}
By Theorem~\ref{thm:Tappflat}, a $g$-zero plane lies in a complete totally geodesic flat immersion $\iota:\R^2\to M$. The conformal sectional-curvature formula~\cite[Chapter 1]{Besse}, applied to a $g$-orthonormal pair tangent to this flat, gives
\begin{equation}\label{eq:conformalflat}
 0\le e^{2f}\secc_{\widehat g}
 =-\Delta_{\R^2}(f\circ\iota)-|\nabla^\perp f|^2.
\end{equation}
The function $f\circ\iota$ is bounded and superharmonic. Parabolicity of $\R^2$, equivalently its Liouville property for positive superharmonic functions~\cite[Theorem 5.1]{Grigoryan}, makes it constant. Equation~\eqref{eq:conformalflat} then forces both the normal-gradient term and the sectional curvature of the plane to vanish. This holds for every original zero plane. The last assertion follows from \eqref{eq:oldall}.
\end{proof}
The same proposition applies to the nested metrics in dimensions twelve and twenty-four. In particular their surviving codimension-one zero sets cannot be removed by a conformal change that retains nonnegative sectional curvature. This is a restriction on the deformation, not an additional topological obstruction.

\section{Fundamental-group orders beyond the flag covers}\label{sec:orders}
The preceding examples separate quasipositivity and almost positivity from the restrictions in Synge's theorem. It is natural to compare their fundamental groups with those allowed by general nonnegative-curvature topology.

\subsection{An order bound when Euler characteristic is nonzero}
We use two classical theorems. First, the Cheeger--Gromoll splitting theorem implies that a closed manifold with nonnegative Ricci curvature and positive-definite Ricci curvature at one point has finite fundamental group~\cite{CheegerGromoll}. Indeed an infinite deck group over a compact base gives a line in the universal cover; its splitting direction would have zero Ricci curvature at every point. Quasipositive sectional curvature satisfies the stated Ricci hypotheses. Second, Gromov's Betti-number theorem~\cite[Theorem 0.2.A]{Gromov} supplies, for each $n$, a constant $B_n$ such that every closed nonnegatively curved $n$-manifold satisfies $\sum_i b_i(M;\Qrat)\le B_n$.

\begin{proposition}\label{prop:Eulerorder}
If $M^n$ is closed and quasipositively curved, then
\begin{equation}\label{eq:Eulerorder}
 |\pi_1(M)|\,|\chi(M)|\le B_n.
\end{equation}
If $M$ is orientable, $n\equiv2\pmod4$, and $\chi(M)\ne0$, then $|\pi_1(M)|\le B_n/2$. In any fixed dimension an unbounded sequence of fundamental-group orders must eventually have zero Euler characteristic, both on the manifolds and on their universal covers.
\end{proposition}
\begin{proof}
The universal Riemannian cover $\widetilde M$ is compact by finiteness of the fundamental group and remains nonnegatively curved. Therefore
\[
 |\pi_1(M)|\,|\chi(M)|
 =|\chi(\widetilde M)|
 \le\sum_i b_i(\widetilde M;\Qrat)\le B_n.
\]
For an oriented closed manifold of dimension $4k+2$, the middle-dimensional rational intersection form is skew-symmetric and nondegenerate, so its rank is even. Pairing the other Betti numbers by Poincar\'e duality shows that $\chi(M)$ is even. If nonzero, its absolute value is at least two. Finally, once $|\pi_1(M)|>B_n$, the integer $\chi(M)$ must vanish, and so does the Euler characteristic of the cover.
\end{proof}

\begin{corollary}\label{cor:fourdim}
For closed quasipositively curved four-manifolds, $|\pi_1(M)|\le B_4$, and $|\pi_1(M)|\le B_4/2$ in the orientable case. In particular no such four-dimensional family has unbounded group order.
\end{corollary}
\begin{proof}
Finite fundamental group gives $b_1(M)=0$. If $M$ is orientable, duality gives $b_3(M)=0$ and $\chi(M)=2+b_2(M)\ge2$. If $M$ is nonorientable, its orientation double cover has Euler characteristic at least two, hence $\chi(M)\ge1$. Apply \eqref{eq:Eulerorder}.
\end{proof}

\subsection{Six-dimensional biquotient universal covers}
DeVito's classification~\cite[Proposition 2.6, Table 1, and Corollary 3.5]{DeVitoBiquotients} has the following consequence: a simply connected compact six-dimensional biquotient is either diffeomorphic to $S^3\times S^3$, or has no odd rational cohomology and Euler characteristic in $\{2,4,6,8\}$. In the latter case $b_2\le3$; for the zero-Euler-characteristic rational homotopy type, Corollary 3.5 identifies the smooth manifold itself with $S^3\times S^3$. We use this classification as an input rather than repeat it.

\begin{proposition}\label{prop:biquotientorder}
Let $M^6$ be closed and quasipositively curved. Suppose its universal cover is diffeomorphic to a compact biquotient other than $S^3\times S^3$. Then
\[
 |\pi_1(M)|\le8,
 \qquad
 |\pi_1(M)|\le4\quad\text{if }M\text{ is orientable}.
\]
Consequently an unbounded-order family with biquotient universal covers must, after finitely many terms, have universal cover diffeomorphic to $S^3\times S^3$.
\end{proposition}
\begin{proof}
Euler multiplicativity and the quoted classification give
$|\pi_1(M)|\chi(M)=\chi(\widetilde M)\in\{2,4,6,8\}$, so the order is at most eight. Rational transfer eliminates the odd cohomology of $M$. If $M$ is orientable, then $\chi(M)=2+2b_2(M)\ge2$, giving the second bound. The last assertion follows by applying this bound to every term whose cover is not $S^3\times S^3$.
\end{proof}
This proposition restricts the smooth universal cover, not the symmetry of the chosen metric. It neither constructs nor excludes a quasipositive metric on a finite quotient of $S^3\times S^3$.

\subsection{Comparison with other quasipositive finite quotients}
Kerr and Tapp construct a quasipositive metric on $T^1S^6$ with a free isometric circle action~\cite[Corollary 1.4]{KerrTapp}. The total space is simply connected, as follows from $S^5\to T^1S^6\to S^6$. For every $m\ge1$, its quotient by the subgroup $C_m$ is therefore an eleven-dimensional quasipositively curved manifold with fundamental group $C_m$. These examples show that quasipositivity alone permits unbounded orders in a fixed dimension. They are orientable and odd-dimensional; Synge's theorem does not rule out positive curvature on them. They should not be identified with examples known to admit no positive metric.

The distinction from Wilking's real projective tangent bundles is complementary. For $r\ge1$,
\[
 P_{\R}T\RP^{2r+1}
\]
has almost-positive curvature~\cite{Wilking,Tapp}, dimension $4r+1$, and is nonorientable. Nonorientability follows by restricting its tangent bundle to the even-dimensional projective fiber: an orientation of the total space would orient that fiber. Synge's odd-dimensional conclusion therefore excludes positive curvature. Its universal cover is the simply connected Stiefel manifold $V_2(\R^{2r+2})=T^1S^{2r+1}$, and the deck transformations independently reverse the two frame vectors. Hence its fundamental group is $C_2\times C_2$, of constant order four. These established families yield unbounded dimensions, whereas \eqref{eq:Eulerorder} and Theorem~\ref{thm:degrees} concern group orders.

\section{Further questions}\label{sec:questions}
The six-dimensional cyclic quotient remains the natural unresolved member of the almost-positive picture. Theorem~\ref{thm:complexrestriction} excludes the entire nested connected-subgroup construction used here, not arbitrary metrics. A successful construction on $X_2$ would have to avoid that restricted flat-plane persistence. Almost positivity would then descend to $Y_2$ only if the additional involution were isometric; it is not legitimate to obtain such invariance by averaging metrics without a curvature argument.

The codimension-one sets $\Psi=1/2$ also distinguish the present constructions from stronger localization questions. Their projection-rank description shows that no additional nested deformation inside the first block subgroup removes them. Proposition~\ref{prop:conformal} excludes conformal removal within nonnegative curvature. Neither fact excludes a different, nonconformal metric with a smaller zero locus on the same manifold.

Finally, the sharp orders three and six depend on the universal cover's Euler characteristic, not on the details of the Cheeger metric. To obtain unbounded orders in a fixed dimension while excluding every positive metric, one must leave these flag covers; by Proposition~\ref{prop:Eulerorder}, a sufficiently large-order family must have Euler characteristic zero. The construction and topology of such a family are separate problems from the finite quotients settled here.

\appendix
\section{Algebraic verification and reproducibility}\label{app:verification}
All existence, completeness, and Hausdorff-dimension assertions are proved in the text. Supplementary exact-arithmetic scripts check the displayed complex and quaternionic identities, normalizer and orientation calculations, the Albert-algebra derivation model, the torus weights and centralizer dimension, the Fourier identities, and the recurrence in Lemma~\ref{e:lem:flagbound}. They also check representative formulas used in the parameter counts. These finite checks are not a replacement for the all-plane classification or the stratified dimension arguments.

The distributed source includes the exact quaternionic and exceptional verification scripts together with their execution results. No numerical curvature samples are used in any proof. The exceptional model is generated from octonion structure constants as $27\times27$ derivation matrices; its commutation calculations are carried out in exact integer or rational arithmetic. The manuscript itself requires no supplementary code to follow the proofs.


\begin{thebibliography}{99}

\bibitem{Baez}
J. C. Baez, \emph{The octonions}, Bull. Amer. Math. Soc. (N.S.) \textbf{39} (2002), 145--205; erratum, \textbf{42} (2005), 213. \doi{10.1090/S0273-0979-01-00934-X}.

\bibitem{BerardBergery}
L. B\'erard-Bergery, \emph{Certaines fibrations d'espaces homog\`enes riemanniens}, Compositio Math. \textbf{30} (1975), 43--61.

\bibitem{Besse}
A. L. Besse, \emph{Einstein Manifolds}, Ergebnisse der Mathematik und ihrer Grenzgebiete (3), vol. 10, Springer-Verlag, Berlin, 1987. \doi{10.1007/978-3-540-74311-8}.

\bibitem{BCR}
J. Bochnak, M. Coste, and M.-F. Roy, \emph{Real Algebraic Geometry}, Ergebnisse der Mathematik und ihrer Grenzgebiete (3), vol. 36, Springer-Verlag, Berlin, 1998. \doi{10.1007/978-3-662-03718-8}.

\bibitem{Cheeger}
J. Cheeger, \emph{Some examples of manifolds of nonnegative curvature}, J. Differential Geom. \textbf{8} (1973), 623--628. \doi{10.4310/jdg/1214431964}.

\bibitem{CheegerGromoll}
J. Cheeger and D. Gromoll, \emph{The splitting theorem for manifolds of nonnegative Ricci curvature}, J. Differential Geom. \textbf{6} (1971), 119--128. \doi{10.4310/jdg/1214430220}.

\bibitem{DeVitoBiquotients}
J. DeVito, \emph{The classification of compact simply connected biquotients in dimension $6$ and $7$}, Math. Ann. \textbf{368} (2017), 1493--1541. \doi{10.1007/s00208-016-1460-8}.

\bibitem{DeVitoJohnson}
J. DeVito and P. Johnson, \emph{Curvature on Eschenburg spaces}, Transform. Groups \textbf{31} (2026), 465--492. \arxiv{2306.15510}.

\bibitem{DeVitoNance}
J. DeVito and E. Nance, \emph{Almost positive curvature on an irreducible compact rank $2$ symmetric space}, Int. Math. Res. Not. IMRN \textbf{2020}, no. 5, 1346--1365. \doi{10.1093/imrn/rny049}.

\bibitem{DeVitoWest}
J. DeVito and J. West, \emph{Almost positively curved generalized Eschenburg spaces}, preprint (2025), \arxiv{2511.17751}.

\bibitem{EschenburgKerin}
J.-H. Eschenburg and M. Kerin, \emph{Almost positive curvature on the Gromoll--Meyer sphere}, Proc. Amer. Math. Soc. \textbf{136} (2008), 3263--3270. \doi{10.1090/S0002-9939-08-09429-X}.

\bibitem{Grigoryan}
A. Grigor'yan, \emph{Analytic and geometric background of recurrence and non-explosion of the Brownian motion on Riemannian manifolds}, Bull. Amer. Math. Soc. (N.S.) \textbf{36} (1999), 135--249. \doi{10.1090/S0273-0979-99-00776-4}.

\bibitem{GromollMeyer}
D. Gromoll and W. Meyer, \emph{An exotic sphere with nonnegative sectional curvature}, Ann. of Math. (2) \textbf{100} (1974), 401--406. \doi{10.2307/1971078}.

\bibitem{Gromov}
M. Gromov, \emph{Curvature, diameter and Betti numbers}, Comment. Math. Helv. \textbf{56} (1981), 179--195. \doi{10.1007/BF02566208}.

\bibitem{Hatcher}
A. Hatcher, \emph{Algebraic Topology}, Cambridge University Press, Cambridge, 2002.

\bibitem{Kerin}
M. Kerin, \emph{Some new examples with almost positive curvature}, Geom. Topol. \textbf{15} (2011), 217--260. \doi{10.2140/gt.2011.15.217}.

\bibitem{KerrTapp}
M. M. Kerr and K. Tapp, \emph{A note on quasi-positive curvature conditions}, Differential Geom. Appl. \textbf{34} (2014), 63--79. \doi{10.1016/j.difgeo.2014.03.009}.

\bibitem{Knapp}
A. W. Knapp, \emph{Lie Groups Beyond an Introduction}, second ed., Progress in Mathematics, vol. 140, Birkh\"auser, Boston, 2002.

\bibitem{Kollross}
A. Kollross, \emph{Octonions, triality, the exceptional Lie algebra $F_4$, and polar actions on the Cayley hyperbolic plane}, Int. J. Math. \textbf{31} (2020), no. 7, 2050051, 28 pp. \doi{10.1142/S0129167X20500512}.

\bibitem{MaOhnita}
H. Ma and Y. Ohnita, \emph{On Lagrangian submanifolds in complex hyperquadrics and isoparametric hypersurfaces in spheres}, Math. Z. \textbf{261} (2009), 749--785. \doi{10.1007/s00209-008-0350-5}.

\bibitem{ONeill}
B. O'Neill, \emph{The fundamental equations of a submersion}, Michigan Math. J. \textbf{13} (1966), 459--469. \doi{10.1307/mmj/1028999604}.

\bibitem{PetersenWilhelm}
P. Petersen and F. Wilhelm, \emph{Examples of Riemannian manifolds with positive curvature almost everywhere}, Geom. Topol. \textbf{3} (1999), 331--367. \doi{10.2140/gt.1999.3.331}.

\bibitem{Synge}
J. L. Synge, \emph{On the connectivity of spaces of positive curvature}, Quart. J. Math. Oxford Ser. \textbf{7} (1936), 316--320. \doi{10.1093/qmath/os-7.1.316}.

\bibitem{Tapp}
K. Tapp, \emph{Quasi-positive curvature on homogeneous bundles}, J. Differential Geom. \textbf{65} (2003), 273--287. \doi{10.4310/jdg/1090511688}.

\bibitem{TappFlats}
K. Tapp, \emph{Flats in Riemannian submersions from Lie groups}, Asian J. Math. \textbf{13} (2009), 459--464. \doi{10.4310/AJM.2009.v13.n4.a2}.

\bibitem{Wallach}
N. R. Wallach, \emph{Compact homogeneous Riemannian manifolds with strictly positive curvature}, Ann. of Math. (2) \textbf{96} (1972), 277--295. \doi{10.2307/1970789}.

\bibitem{Wilhelm}
F. Wilhelm, \emph{An exotic sphere with positive curvature almost everywhere}, J. Geom. Anal. \textbf{11} (2001), 519--560.

\bibitem{Wilking}
B. Wilking, \emph{Manifolds with positive sectional curvature almost everywhere}, Invent. Math. \textbf{148} (2002), 117--141. \doi{10.1007/s002220100190}.

\end{thebibliography}
\end{document}